\documentclass{svjour3}     
\usepackage{graphicx}
\usepackage{amsmath}
\usepackage{epstopdf} 
\usepackage[ansinew]{inputenc}
\usepackage{amssymb}
\usepackage{makeidx}
\usepackage{float}
\restylefloat{table}

\title{Rodas6P and Tsit5DA - two new Rosenbrock-type methods for DAEs}
\author{Gerd Steinebach}

\institute{Gerd Steinebach \at
              Bonn-Rhein-Sieg University of Applied Sciences, Germany\\
              \email{Gerd.Steinebach@h-brs.de}           
}

\begin{document}
\maketitle

\begin{abstract}
Two new Rosenbrock methods for solving index-1 differential algebraic equations are presented. 
{\tt Rodas6P} is a sixth-order method based on the same design principles as the {\tt Rodas3P}, {\tt Rodas4P}, and {\tt Rodas5P} methods, \cite{rodas4p2,rodas5p,steijulia}. 
{\tt Tsit5DA} is based on an explicit solution approach for the differential equations and a linear-implicit approach for the algebraic equations. 
Such a fourth-order method has already been presented in \cite{renrostei}.  {\tt Tsit5DA} now provides a significantly improved fifth-order method which is based on {\tt Tsit5} \cite{tsit5}.
The theoretical properties of the new methods are verified by some order tests and benchmarks.
\end{abstract}

\section{Rosenbrock-type methods}

Roche \cite{roche} first considered Rosenbrock-Wanner (ROW) methods for differential algebraic equations (DAEs) of type
\begin{eqnarray}
 y' &=& f(t,y,z) \, , \; y(t_0)=y_0 \, , \label{eq:dae1} \\
 0 &=& g(t,y,z) \, , \; z(t_0)=z_0 \, . \label{eq:dae2} 
\end{eqnarray}
Here, consistent intial values with $g(t_0,y_0,z_0)=0$ are assumed and the index-1 assumption of a regular matrix $\frac{\partial g}{\partial z}$
in the neighborhood of the soluton $y(t)$, $z(t)$ must be fullfilled. We assume that dimensions of vector $y$ is $n_f$ and of vector $z$ is $n_g$.
Hairer \& Wanner \cite{hairer} formulated the methods for non-autonomous problems in so-called mass-matrix form of dimension $n_f+n_g$:
\begin{equation}
M \, y' = f(t,y) \, , \; y(t_0)=y_0 \, . \label{eq:dae} 
\end{equation}
If the mass matrix $M$ is singular, it is a DAE problem, and with the help of a singular value decomposition of $M$, 
equation (\ref{eq:dae}) can be reformulated as (\ref{eq:dae1},\ref{eq:dae2}).

A ROW scheme with stage-number $s$ for problem (\ref{eq:dae}) is defined by:
\begin{eqnarray}
(M - h \, \gamma \, f_y) k_i &=& h f(t_0+ \alpha_i h, y_0 + \sum_{j=1}^{i-1} \alpha_{ij} k_j) \nonumber \\
&& + h \, f_y \sum_{j=1}^{i-1} \gamma_{ij}  k_j + h^2 \gamma_i f_t \; ,  \qquad i=1,...,s,  \label{eq_row1} \\
y_1 &=& y_0 + \sum_{i=1}^{s} b_i k_i,  \label{eq_row2} \\
\mbox{with} \quad f_y &=& \frac{\partial f}{\partial y}(t_0,y_0) \; , \quad f_t = \frac{\partial f}{\partial t}(t_0,y_0) \, . \nonumber
\end{eqnarray}
$h$ is the stepsize and $y_1$ is the approximation to the solution $y(t_0+h)$.
The coefficients of the method are $\gamma$, $\alpha_{ij}$, $\gamma_{ij}$, and $b_i$ define the weights. Moreover, it holds $\alpha_i = \sum_{j=1}^{i-1} \alpha_{ij}$ and
$\gamma_i = \gamma + \sum_{j=1}^{i-1} \gamma_{ij}$.\\
The Jacobian matrix $f_y$ must be reevaluated at each time step. Through the use of modern automatic differentiation such as 
{\tt ForwarDiff.jl}, this disadvantage of ROW methods is becoming increasingly less significant.
W methods avoid the need for an up-to-date Jacobian matrix at each time step. However, in order to derive higher-order W methods, an immensely 
higher number of order conditions must be fulfilled compared to ROW methods. Jax \cite{jax2}, for example, requires a stage number of $s=7$ for a third-order W method. 
Efficient third-order ROW methods like {\tt ROS3P} \cite{lang2} however, can already be achieved  with $s=3$ stages.

The linear systems that arise in equation (\ref{eq_row1}) are usually solved by performing an LU decomposition of the matrix $E = M - h \, \gamma \, f_y$ at each time step, 
so that the subsequent $s$ back substitutions can be computed efficiently. Even though the sparse nature of the matrix E is usually exploited in large systems, 
LU decomposition can be computationally expensive.

A compromise between full W methods and ROW schemes is the approach presented in \cite{renrostei} for problems of type (\ref{eq:dae1},\ref{eq:dae2}):
\begin{eqnarray}
l_i &=& h f(t_0+ \alpha_i h, y_0 + \sum_{j=1}^{i-1} \alpha_{ij} l_j,z_0 + \sum_{j=1}^{i-1} \alpha_{ij} k_j) \label{eq_rowda1} \\
- \gamma \, g_z k_i &=& g(t_0+ \alpha_i h, y_0 + \sum_{j=1}^{i-1} \alpha_{ij} l_j,z_0 + \sum_{j=1}^{i-1} \alpha_{ij} k_j) \nonumber \\
&& + g_y \sum_{j=1}^{i} \gamma_{ij}  l_j + h \gamma_i g_t + g_z \sum_{j=1}^{i-1} \gamma_{ij}  k_j \; ,  \qquad i=1,...,s,  \label{eq_rowda2} \\
\begin{pmatrix} y_1\\z_1 \end{pmatrix} &=& y_0 + \sum_{i=1}^{s} b_i \begin{pmatrix} l_i\\k_i\end{pmatrix},  \label{eq_rowda3} \\
\mbox{with} \quad g_y &=& \frac{\partial g}{\partial y}(t_0,y_0,z_0) \; , \quad g_z = \frac{\partial g}{\partial z}(t_0,y_0,z_0) \; , \quad 
g_t = \frac{\partial g}{\partial t}(t_0,y_0,z_0) \, . \nonumber
\end{eqnarray}
For pure differential equations $y'=f(t,y)$, this approach is equivalent to an explicit Runge-Kutta method.
With this approach, only the LU decomposition of the matrix $g_z$ with dimension $n_g$ needs to be calculated.
A somewhat more general approach, which also includes approximations to the matrices $f_y$ and $f_z$, is considered in \cite{jaxstei}. 
However, this again would require the LU decomposition of a matrix of dimension $n_f+n_g$.

Lang \cite{lang} provides an excellent overview of the development and current status of Rosenbrock methods. 
There already exist a large number of methods for approach (\ref{eq_row1},\ref{eq_row2}).
The Rosenbrock methods currently implemented in the Julia package {\tt OrdinaryDiffEq.jl} \cite{julia} are described in \cite{steijulia}.
The fifth-order method {\tt Rodas5P} \cite{rodas5p} is one of the schemes recommended for solving stiff differential equations and differential algebraic equations with medium to high accuracy requirements.
The question arises whether it is possible to develop a method of order $p=6$ with properties comparable to those of {\tt Rodas5P}.

For non-stiff problems, the explicit Runge-Kutta method {\tt Tsit5} \cite{tsit5} is recommended in {\tt OrdinaryDiffEq.jl}, among others. 
One way to solve DAEs using explicit methods is to solve the algebraic equations each time the right-hand side is called.
To reduce this effort, the approach (\ref{eq_rowda1},\ref{eq_rowda2},\ref{eq_rowda3}) can be appplied.
In \cite{renrostei}  method {\tt RKF4DA} of order $p=4$ was presented , which is based on the explicit Runge-Kutta-Fehlberg method {\tt RKF45}, \cite{fehlberg}.
This raises the question of whether it is possible to derive an analogous method of order $p=5$ based on the more efficient explicit Runge-Kutta method {\tt Tsit5}, \cite{tsit5}.

\section{Construction of Rodas6P}
Like the original {\tt Rodas} method \cite{hairer}, {\tt Rodas6P} should be stiffly accurate and it's embedded scheme of order $\hat p=5$ with weigths 
$\hat b_i$, $i=1,\ldots,s-1$, too:
\begin{eqnarray}
b_i&=&\beta_{s, i} \quad \mbox{for } i=1,\ldots,s-1 , \quad b_s=\gamma ,\quad \alpha_s=1 \; .\label{stiffly}\\
\hat b_i&=&\beta_{s-1, i} \quad \mbox{for } i=1,\ldots,s-2 , \quad \hat b_{s-1}=\gamma ,
\quad \alpha_{s-1}=1 \ . \label{stiffly1}
\end{eqnarray}
The coefficients $\beta$ are given by $\beta_{i,j}=\alpha_{i,j}+\gamma_{i,j}$, $\beta_{i,i}=\gamma$.
For a ROW method of type (\ref{eq_row1},\ref{eq_row2}) with order $p=6$, 130 order conditions must be satisfied. 
These are listed together with their Butcher trees in Table \ref{tab:rodas6p} in the appendix. Note, that coefficients 
$w_{i,j}$ arise from matrix $ W = \mathcal{B}^{-1}$ with matrix $\mathcal{B}$ consisting of coefficients $\beta_{i,j}$.
For the embedded method with weights $\hat b_i$ the first 40 conditions of Table \ref{tab:rodas6p} muss be fullfilled.

Moreover, {\tt Rodas6P} should have similar properties to {\tt Rodas5P}. The letter {\tt P} stands for reduced order reduction in the context of
the Prothero-Robinson model and semi-discretized parabolic partial differential equations, \cite{oro,rang2,scholz}. Therefore, the same additional 9 conditions as given in
Table 2.3 in \cite{rodas5p} are required. They also include improvements for the method when applied to index-2 problems or when the Jacobian $f_y$
is not exact.

In order to fulfill all 179 conditions, a stage number of $s=16$ was required. 
The coefficients could be calculated using the Levenberg-Marquardt and Trust-Region methods of the package  {\tt NonlinearSolve.jl}.
By choosing $\gamma=0.26$, A-stability of the 6th-order method and the embedded 5th-order method could be achieved. Due to the stiffly accurate property, the
methods are L-stable as well.

To equip {\tt Rodas6P} with internal interpolation of order five for dense output, analogous to Hairer \& Wanner \cite{hairer}, equations (\ref{eq:interp1}), (\ref{eq:interp2})
are applied:
\begin{eqnarray}
 y(t_0 + \tau \, h) &=& y_0 + \sum_{i=1}^s b_i(\tau) k_i \; , \quad \tau \in [0,1] \, ,\label{eq:interp1} \\
 b_i(\tau) &=& \tau (b_i-c_i) + \tau^2 (c_i-d_i) + \tau^3 (d_i-e_i) + \tau^4 (e_i-f_i) + \tau^5 f_i \, .\label{eq:interp2}
\end{eqnarray}
The computation of each set of coeffcients $c_i$, $d_i$, $e_i$ and $f_i$ for $i=1,...,s$ requires to solve 40 order conditions.
When coefficients $\alpha_{i,j}$, $\beta_{i,j}$ are given, the resulting 40 equations for e.g. $c_i$, $i=1,..,s$ are linear. But for $s=16$ the system is not solvable.
Therefore, the number of stages was increased to $s=19$. This results in 102 new degrees of freedom for the coefficients 
$\alpha_{i,j}$, $\beta_{i,j}$, $i=17,18,19$, $j=1,...,i-1$
and 76 degrees of freedom for $c_i$, $d_i$, $e_i$, $f_i$, $i=1,...,19$.
In summary, we now have a nonlinear system of 160 equations for 178 unknowns. 
In addition, some of the further conditions reducing order reduction specified by Scholz \cite{scholz}, could also be met for the interpolation according to (\ref{eq:interp1},\ref{eq:interp2}).
An implementation of {\tt Rodas6P} and it's coefficients can be found in {\tt OrdinaryDiffEqRosenbrock.jl}, which is a sublibrary of {\tt OrdinaryDiffEq.jl}.

\section{Construction of Tsit5DA}
{\tt Tsit5} is an explicit Runge-Kutta method of order 5(4) with $s=7$ stages, \cite{tsit5}.
The seventh stage is only necessary for the embedded method, so that $b_7=0$ applies.
Table \ref{tab:tsit5da} in the appendix shows the 63 order conditions up to order $p=5$ for a method of type (\ref{eq_rowda1},\ref{eq_rowda2},\ref{eq_rowda3}).
It is not possible to take all coefficients $\alpha_{i,j}$, $b_i$, $\hat b_i$ from {\tt Tsit5} and simply determine 22 new coefficients
$\beta_{i,j}$, $i=1,...,7$, $j=1,...,i-1$ and $\beta_{i,i}=\gamma$ to fullfill all order conditions.
We therefore choose $s=12$ and set $\alpha_{i,j}$ as follows:
\begin{equation}
\alpha = \left( \begin{array}{cccccccccccc}
                        0&0&0&0&0&0&0&0&0&0&0&0 \\
                        2 \gamma &0&0&0&0&0&0&0&0&0&0&0 \\
                        a &0&0&0&0&0&0&0&0&0&0&0 \\
                        \alpha_{21} &0&0&0&0&0&0&0&0&0&0&0 \\
                        \alpha_{31} &0&0&\alpha_{32}&0&0&0&0&0&0&0&0 \\
                        \alpha_{41} &0&0&\alpha_{42}&\alpha_{43}&0&0&0&0&0&0&0 \\
                        \alpha_{51} &0&0&\alpha_{52}&\alpha_{53}&\alpha_{54}&0&0&0&0&0&0 \\
                        \alpha_{61} &0&0&\alpha_{62}&\alpha_{63}&\alpha_{64}&\alpha_{65}&0&0&0&0&0 \\
                        b_1 &0&0&b_2&b_3&b_4&b_5&b_6&0&0&0&0 \\
                        \hat b_1 &0&0&\hat b_2&\hat b_3&\hat b_4&\hat b_5&\hat b_6&\hat b_7&0&0&0 \\
                        b_1 &0&0&b_2&b_3&b_4&b_5&b_6&0&0&0&0 \\
                        \hat b_1 &0&0&\hat b_2&\hat b_3&\hat b_4&\hat b_5&\hat b_6&\hat b_7-\gamma&0&\gamma&0 
\end{array} \right)
\end{equation}
Thus we have two free coefficients $a$ and $\gamma$, all other values are taken from {\tt Tsit5}.
Furthermore, we set the new weights of {\tt Tsit5DA} to $\hat b_i = \alpha_{12,i}$ for $i=1,..,12$ and $b_i=\alpha_{11,i}$ for $i=1,..,9$ and $b_{10}=-\gamma$, $b_{11}=0$, $b_{12}=\gamma$.
Since we want {\tt Tsit5DA} to be stiffly accurate, we set again the coefficients $\beta_{s-1,i}$ and $\beta_{s,i}$ according to (\ref{stiffly},\ref{stiffly1}).
This property guarantees convergence, because $|R(\infty)|=0<1$ applies to the stability function.
In this way, we now have 47 free coefficients: $\gamma$, $a$, $\beta_{2,1},...,\beta_{10,9}$.

Note, that the application of this scheme to pure differential equations 
is equivalent to {\tt Tsit5} and that all order conditions of Table \ref{tab:tsit5da} with only thin knots are fullfilled.
To compute the $\beta$ coeffcients we apply the simplifying assumptions
\begin{equation}
\sum_j w_{i,j} \alpha_j^2 = 2 \alpha_i \quad \mbox{for} \; i = 2,...,s \, .\label{eq:simple} 
\end{equation}
This means that all trees with branches of the form $\vcenter{\hbox{\includegraphics[width = 0.055\textwidth]{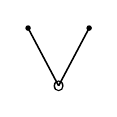}}}$ can be traced back to those with 
branches of the form $\vcenter{\hbox{\includegraphics[width = 0.055\textwidth]{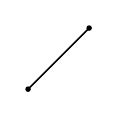}}}$.
We now set $a=0.4$, $\gamma=0.15$ and can calculate all coefficients using again {\tt NonlinerSolve.jl}.

For dense output of order four we apply equations (\ref{eq:interp1},\ref{eq:interp2}) with $f_i=0$. 
Due to the simplifying assumptions (\ref{eq:simple}), the necessary 18 conditions for order four are reduced to 12, and we can calculate 
each coefficient set $c_i$, $d_i$, $e_i$ by solving a linear system. Table \ref{tab:coeff} in the appendix lists all coefficients of {\tt Tsit5DA}.

\section{Benchmarks and Conclusion}
First, we perform two simple order tests taken from \cite{rodas5p}. In test one we treat the DAE
\begin{equation}
y_1' = \frac{y_2}{y_1} \; ,\quad
0 = \frac{y_1}{y_2} - t \, , \label{eq:prob1}
\end{equation}
with exact solution $y_1(t)=\ln(t)$, $y_2(t)=\frac{\ln(t)}{t}$ in the time interval $t \in [2,4]$.
We apply methods {\tt Rodas3P} \cite{steijulia}, {\tt Rodas4P} \cite{rodas4p}, {\tt Rodas5P} \cite{rodas5p}, {\tt RKF4DA} \cite{renrostei}
and the new methods {\tt Rodas6P}, {\tt Tsit5DA} to this problem with different constant stepsizes. Table \ref{tab:prob1} shows the obtained absolute errors at timt $t=4$
and the numerically determined orders of the methods. We can clearly see that the theoretical orders are being achieved.

\begin{table}
\caption{Numerical errors and orders of the various methods for problem (\ref{eq:prob1}).}
\label{tab:prob1}
\begin{center}
\begin{tabular}{|r|r|r|r|r|r|r|}
  \hline
  \textbf{stepsize} & \textbf{Rodas3P} & \textbf{Rodas4P} & \textbf{Rodas5P} & \textbf{Rodas6P} & \textbf{RKF4DA} & \textbf{Tsit5DA} \\
  \hline
  1.25e-01 & 3.18e-05 & 3.10e-07 & 2.93e-08 & 5.03e-10 & 3.02e-06 & 1.51e-07 \\
  6.25e-02 & 4.05e-06 & 1.79e-08 & 8.56e-10 & 7.25e-12 & 1.58e-07 & 4.03e-09 \\
  3.12e-02 & 5.10e-07 & 1.08e-09 & 2.59e-11 & 1.09e-13 & 9.06e-09 & 1.22e-10 \\
  1.56e-02 & 6.41e-08 & 6.64e-11 & 8.01e-13 & 3.77e-15 & 5.43e-10 & 3.79e-12 \\
  7.81e-03 & 8.02e-09 & 4.12e-12 & 2.93e-14 & 4.44e-15 & 3.32e-11 & 1.19e-13 \\
  \hline
  \hline
6.25e-02  & 2.97 & 4.11 & 5.10 & 6.11 & 4.26 & 5.22 \\
3.12e-02  & 2.99 & 4.05 & 5.05 & 6.06 & 4.13 & 5.04 \\
1.56e-02  & 2.99 & 4.02 & 5.02 & 4.85 & 4.06 & 5.01 \\
7.81e-03  & 3.00 & 4.01 & 4.77 & -0.23 & 4.03 & 4.99 \\
  \hline
  \hline
  \textbf{stepsize} & \multicolumn{6}{c|}{embedded schemes} \\
  \hline
  1.25e-01 & 1.05e-04 & 8.09e-06 & 1.13e-06 & 1.62e-08 & 1.78e-05 & 1.99e-03 \\
  6.25e-02 & 2.68e-05 & 8.78e-07 & 6.60e-08 & 4.82e-10 & 2.30e-06 & 4.13e-05 \\
  3.12e-02 & 6.74e-06 & 1.01e-07 & 4.00e-09 & 1.47e-11 & 2.89e-07 & 1.77e-08 \\
  1.56e-02 & 1.69e-06 & 1.22e-08 & 2.46e-10 & 4.55e-13 & 3.62e-08 & 1.38e-09 \\
  7.81e-03 & 4.23e-07 & 1.49e-09 & 1.53e-11 & 1.09e-14 & 4.53e-09 & 9.79e-11 \\
  \hline
  \hline
6.25e-02  & 1.98 & 3.20 & 4.10 & 5.07 & 2.96 & 5.59 \\
3.12e-02  & 1.99 & 3.11 & 4.05 & 5.03 & 2.99 & 11.19 \\
1.56e-02  & 2.00 & 3.06 & 4.02 & 5.02 & 3.00 & 3.68 \\
7.81e-03  & 2.00 & 3.03 & 4.01 & 5.39 & 3.00 & 3.82 \\
  \hline
\end{tabular}
\end{center}
\end{table}

The second test is the Prothero-Robinson model 
\begin{equation} y' = -\lambda (y - g) + g'\label{eq:prob2} \end{equation}
with function $g(t)= 10 - (10+t) e^{-t}$ and
stiffness parameter $\lambda = 10$ in the time interval $t \in [0,2]$. The solution is $y(t)=g(t)$ and Table \ref{tab:prob2} shows the 
obtained results  analogous to Table \ref{tab:prob1}.
Here, too, we see that the theoretical orders are achieved and we see that the step sizes for {\tt RKF4DA} and {\tt Tsit5DA} 
must be small enough to lie within the stability range of the methods.

\begin{table}
\caption{Numerical errors and orders of the various methods for problem (\ref{eq:prob2}).}
\label{tab:prob2}
\begin{center}
\begin{tabular}{|r|r|r|r|r|r|r|}
  \hline
  \textbf{stepsize} & \textbf{Rodas3P} & \textbf{Rodas4P} & \textbf{Rodas5P} & \textbf{Rodas6P} & \textbf{RKF4DA} & \textbf{Tsit5DA} \\
  \hline
  5.00e-01 & 8.89e-03 & 6.31e-05 & 1.93e-05 & 9.95e-07 & 5.02e+01 & 8.44e+02 \\
  2.50e-01 & 1.28e-03 & 4.31e-06 & 8.65e-07 & 1.36e-08 & 2.19e-03 & 1.81e-03 \\
  1.25e-01 & 1.80e-04 & 2.87e-07 & 2.92e-08 & 9.71e-11 & 5.43e-05 & 1.63e-05 \\
  6.25e-02 & 2.46e-05 & 1.85e-08 & 8.66e-10 & 5.76e-13 & 1.34e-06 & 2.30e-07 \\
  3.12e-02 & 3.25e-06 & 1.18e-09 & 2.49e-11 & 3.73e-14 & 3.68e-08 & 4.19e-09 \\
  1.56e-02 & 4.21e-07 & 7.43e-11 & 7.25e-13 & 5.33e-15 & 1.07e-09 & 9.26e-11 \\
  7.81e-03 & 5.36e-08 & 4.67e-12 & 2.49e-14 & 0.00e+00 & 3.24e-11 & 2.35e-12 \\
  \hline
  \hline
  2.50e-01 & 2.80 & 3.87 & 4.48 & 6.20 & 14.48 & 18.83 \\
  1.25e-01 & 2.83 & 3.91 & 4.89 & 7.13 & 5.33 & 6.80 \\
  6.25e-02 & 2.87 & 3.95 & 5.07 & 7.40 & 5.34 & 6.14 \\
  3.12e-02 & 2.92 & 3.98 & 5.12 & 3.95 & 5.19 & 5.78 \\
  1.56e-02 & 2.95 & 3.99 & 5.10 & 2.81 & 5.10 & 5.50 \\
  7.81e-03 & 2.97 & 3.99 & 4.87 & Inf & 5.05 & 5.30 \\
  \hline
  \hline
  \textbf{stepsize} & \multicolumn{6}{c|}{embedded schemes} \\
  \hline
  5.00e-01 & 1.74e-03 & 1.17e-04 & 1.26e-04 & 2.43e-04 & 7.31e+03 & 3.98e+01 \\
  2.50e-01 & 3.87e-04 & 5.05e-06 & 6.08e-06 & 1.44e-05 & 5.02e-03 & 1.61e-04 \\
  1.25e-01 & 8.86e-05 & 4.23e-08 & 2.82e-07 & 6.69e-07 & 1.80e-04 & 1.54e-05 \\
  6.25e-02 & 2.09e-05 & 5.42e-08 & 1.33e-08 & 2.63e-08 & 6.12e-06 & 8.87e-07 \\
  3.12e-02 & 5.04e-06 & 1.02e-08 & 6.58e-10 & 9.27e-10 & 2.64e-07 & 4.75e-08 \\
  1.56e-02 & 1.24e-06 & 1.52e-09 & 3.50e-11 & 3.08e-11 & 1.34e-08 & 2.67e-09 \\
  7.81e-03 & 3.06e-07 & 2.05e-10 & 1.98e-12 & 9.84e-13 & 7.49e-10 & 1.57e-10 \\
  \hline
  \hline
  2.50e-01 & 2.17 & 4.53 & 4.38 & 4.08 & 20.47 & 17.92 \\
  1.25e-01 & 2.13 & 6.90 & 4.43 & 4.42 & 4.80 & 3.39 \\
  6.25e-02 & 2.08 & -0.36 & 4.41 & 4.67 & 4.88 & 4.11 \\
  3.12e-02 & 2.05 & 2.40 & 4.33 & 4.83 & 4.54 & 4.22 \\
  1.56e-02 & 2.03 & 2.76 & 4.23 & 4.91 & 4.30 & 4.15 \\
  7.81e-03 & 2.02 & 2.89 & 4.14 & 4.97 & 4.16 & 4.09 \\
  \hline
\end{tabular}
\end{center}
\end{table}

\begin{figure}[t]
 \centering
 \includegraphics[width=1.0\textwidth]{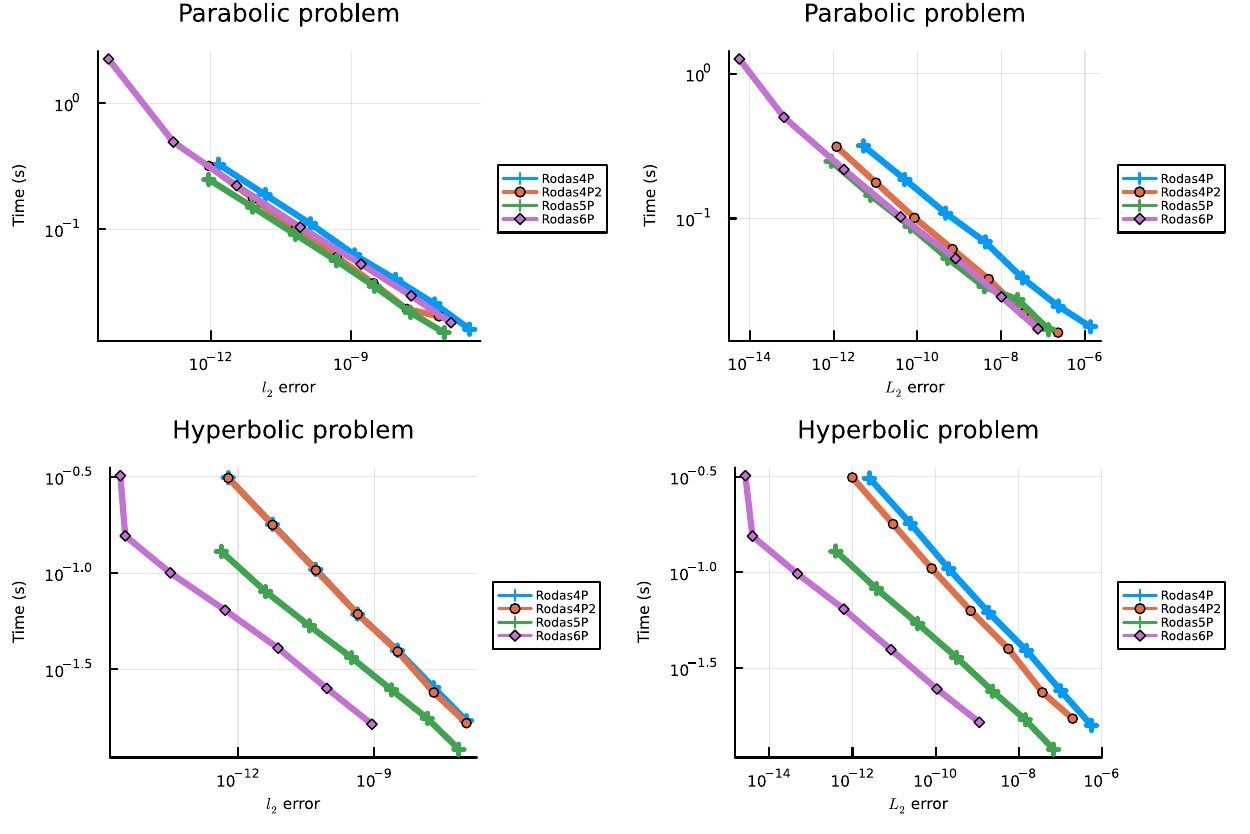}
 \caption{Work-precision diagrams for parabolic and hyperbolic problems.}\label{fig:bench1}
\end{figure}

Next, we show the potential of {\tt Rodas6P}, which is particularly apparent in semi-discretized partial differential equations. 
We again use two benchmarks from \cite{rodas5p} for this purpose.
The parabolic problem reads
\begin{equation}
\frac{\partial u}{\partial t}= \frac{\partial^2 u}{\partial x^2} +u^2 +h(x,t) \, , \;  x\in[-1,1] \, , \; t \in [0,1]
\label{eq:parabol} 
\end{equation}
where function $h(x,t)$ is chosen in order to get the solution $u(x,t)= x^3 \cdot e^t $. 
Initial values and Dirichlet boundary condition are taken from this solution and
discretization in space by $\frac{\partial^2}{\partial x^2}u(x_i,t)=\frac{u(x_{i-1},t)-2 u(x_i,t)+u(x_{i+1},t)}{\Delta x^2}$ leads to a ODE system of
dimension $n_x$.
Figure \ref{fig:bench1} shows the work-precision diagram obtained by the tool {\tt DiffEqDevTools.jl} with $n_x=250$ space discretization points.
The graphs on the left refer to the $l_2$-error taken from the solution at every timestep and on the right refer to the
$L_2$-error taken at 100 evenly spaced points via interpolation, which reflects the error of the dense output formulae.

The hyperbolic problem is given by
\begin{equation}
\frac{\partial u}{\partial t}= -\frac{\partial u}{\partial x} + g(x,t) \, , \;  x\in[0,1] \, , \; t \in [0,1]
\label{eq:hyperbol} 
\end{equation}
where function $g(x,t)$ is chosen in order to get the solution $u(x,t)= \frac{1+x}{1+t}$. Again, $n_x=250$ space points are applied and the
work-precission diagram is also shown in Figure \ref{fig:bench1}.
In addition to the methods mentioned above, {\tt Rodas4P2} \cite{rodas4p2}, which is a minor improvement on {\tt Rodas4P}, is also included.
The results show that {\tt Rodas6P} achieves higher accuracy and that computing time savings are possible, especially for the hyperbolic problem.

It should be noted that there are many examples in which {\tt Rodas5P} works more efficiently than {\tt Rodas6P}. 
The disadvantage of the higher number of stages of $s=19$ in {\tt Rodas6P} compared to $s=8$ in {\tt Rodas5P} is only compensated for if a significantly 
lower number of time steps of the higher-order method is achieved and the dimension of the DAE system is correspondingly high.

The {\tt RKF4DA} and {\tt Tsit5DA} methods were not tested for these examples because the semi-discretized PDEs lead to stiff problems.

\begin{figure}[t]
 \centering
 \includegraphics[width=0.7\textwidth]{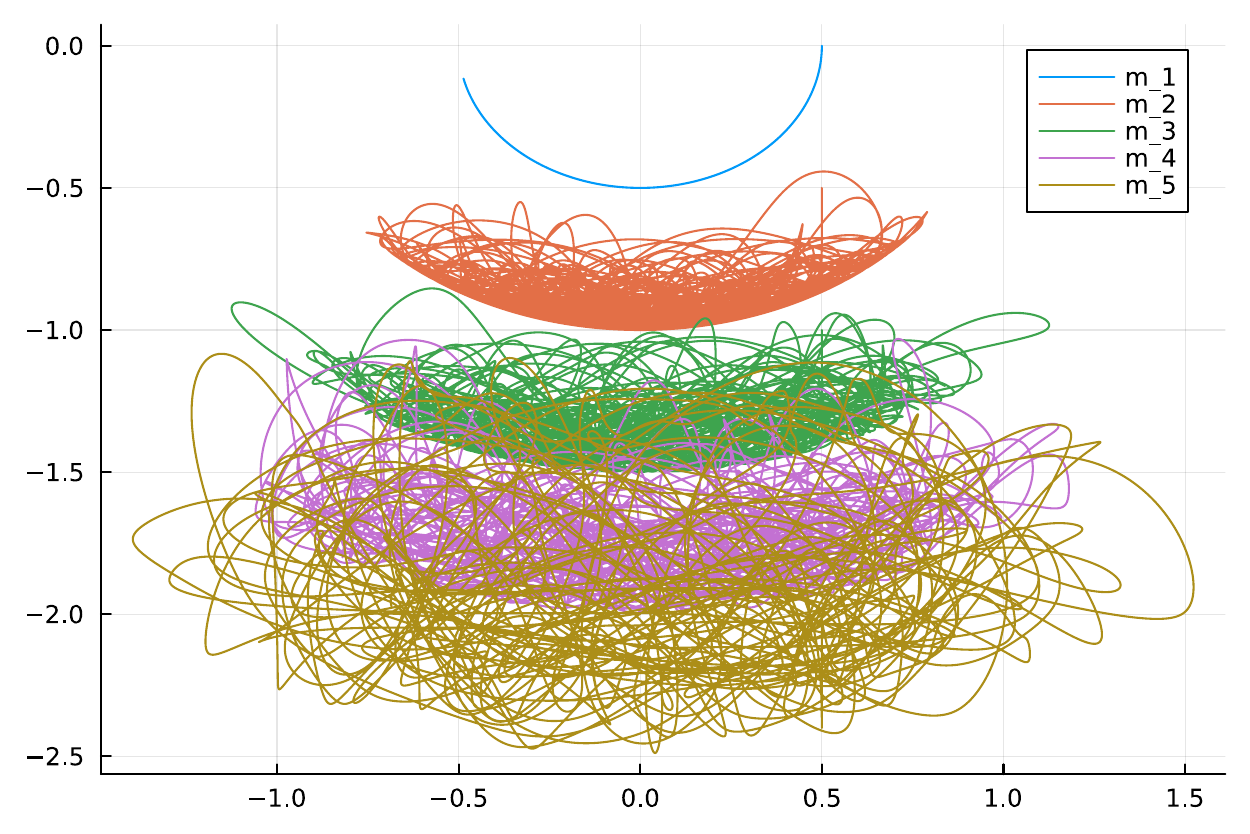}
 \caption{Trajectories of a pendulum with 5 masses.}\label{fig:pendel}
\end{figure}

To illustrate the intended use of {\tt Tsit5DA}, we take the muli-pendulum, implemented as a DAE system, as an example.
It is a dynamic oscillating system without any stiff components.
If we assume that the mass $m_i$ is connected by massless rods of lengths $L_i$ and $L_{i+1}$ with masses $m_{i-1}$ and $m_{i+1}$, we can describe the motion by:
\begin{eqnarray}
m_i x_i'' &=& \lambda_i (x_i - x_{i-1}) - \lambda_{i+1} (x_{i+1}-x_i) \; , \label{eq:pendel1} \\
m_i y_i'' &=& -g + \lambda_i (y_i - y_{i-1}) - \lambda_{i+1} (y_{i+1}-y_i)\; , \label{eq:pendel2} \\
0 &=& (x_i' - x_{i-1}')^2 + (y_i' - y_{i-1}')^2 \nonumber \\
&& + (x_i - x_{i-1})(x_i'' - x_{i-1}'') + (y_i - y_{i-1})(y_i'' - y_{i-1}'') \; . \label{eq:pendel3}
\end{eqnarray}
Here, $x_i(t)$, $y_i(t)$ describe the Cartesian coordinates of the mass $m_i$, and the constraint $(x_i-x_{i-1})^2+(y_i-y_{i-1})^2-L_i^2=0$ 
has already been differentiated twice with respect to time, yielding an index-1 DAE system.
The equations must still be transformed into a first-order system and the second derivatives form (\ref{eq:pendel1}), (\ref{eq:pendel2}) inserted into (\ref{eq:pendel3}). 
A multi-pendulum with $n$ masses then results in a DAE system of dimension $5n$, of which $n$ are algebraic equations. 
The equations for the first and last masses must still be adjusted accordingly.
As an example, Figure \ref{fig:pendel} shows the trajectories of a multi-pendulum with 5 masses within time $t \in [0,100]$.
The trajectories depend on the choice of initial values.

Since {\tt Tsit5DA} is not yet implemented in {\tt OrdinaryDiffEq.jl}, we use our own implementations of all tested methods and do not use work-precision diagrams
from {\tt DiffEqDevTools.jl}. Table \ref{tab:pendel} shows a summary of the results achieved.
The table lists: number of successful time steps NSUCC, number of rejected time steps NFAIL, number of function evaluations NFCN, errors ERR, and computing time CPU.
The maximum difference between the simulated and exact lengths $\sum L_i$ was chosen as the error measure ERR.
The upper part of the table contains the results for RelTol=AbsTol = 1.0e-7 and the lower part for tolerances of 1.0e-8.
The results show that {\tt Tsit5DA} is very well suited for this type of problem and can be an alternative to established ROW methods.

In summary, it can be said that {\tt Rodas6P} and {\tt Tstit5DA} provide two additional Rosenbrock methods that can be used to further increase efficiency in certain applications.
{\tt Rodas6P} complements the Rodas family of methods and can lead to computing time savings compared to {\tt Rodas5P}, 
especially for large DAE systems that arise, for example, from the semi-discretisation of partial differential equations.
But also the highly dynamic example of the multi-pendulum shows that {\tt Rodas6P} can achieve significant savings in time steps compared to lower order methods.

{\tt Tsit5DA} provides a significant improvement of {\tt RKF4DA}, both are Rosenbrock methods with an explicit approach for the differential part of the DAE system.
Possible areas of application are systems without stiff solution components. The potential of the method was also demonstrated with the example of the multi-pendulum.

\begin{table}
\caption{Numerical for the multi-pendulum problem with 5 masses.}
\label{tab:pendel}
\begin{center}
\begin{tabular}{|r|r|r|r|r|r|r|}
  \hline
   & \textbf{Rodas3P} & \textbf{Rodas4P} & \textbf{Rodas5P} & \textbf{Rodas6P} & \textbf{RKF4DA} & \textbf{Tsit5DA} \\
  \hline
  NSUCC &171832&96966&60527&28771&206807&62205\\
  NFAIL  &30&614&899&1053&211&1204\\
  NFCN  &859280&584866&490509&476131&1241897&696295\\
  ERR  &5.9e-4&1.9e-4&5.7e-5&1.0e-5&8.9e-4&2.2e-4\\
  CPU  &14.1&8.8&5.9&3.6&4.1&2.5\\
  \hline
  NSUCC &379666&176137&109201&45629&464476& 106640\\
  NFAIL &4&414&699&734&50&906 \\
  NFCN &1898346&1058892&878501&741074&2787106& 1182100\\
  ERR &1.5e-5&1.9e-5&2.4e-6&1.0e-6&4.0e-5&7.6e-6 \\
  CPU &32.3&14.9&10.3&5.7&8.5&3.8 \\
  \hline
\end{tabular}
\end{center}
\end{table}

\newpage
\appendix
\section{Appendix}
\begin{table}[H] 
\caption{Order conditions for a ROW mthods up to order $p=6$.}
\label{tab:rodas6p}
\end{table}
\vspace*{-0.5cm}
\begin{tabular}{lcl}
\hline
1. & $\vcenter{\hbox{\includegraphics[width = 0.075\textwidth]{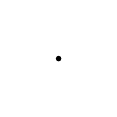}}}$ & $\;\sum b_i= \frac{1}{1}$\\
2. & $\vcenter{\hbox{\includegraphics[width = 0.075\textwidth]{baum6p_2.pdf}}}$ & $\;\sum b_i\beta_{ij}= \frac{1}{2}$\\
3. & $\vcenter{\hbox{\includegraphics[width = 0.075\textwidth]{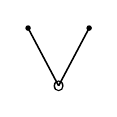}}}$ & $\;\sum b_iw_{ij}\alpha_{jk}\alpha_{jl}= \frac{1}{1}$\\
4. & $\vcenter{\hbox{\includegraphics[width = 0.075\textwidth]{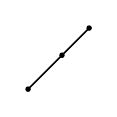}}}$ & $\;\sum b_i\beta_{ij}\beta_{jk}= \frac{1}{6}$\\
5. & $\vcenter{\hbox{\includegraphics[width = 0.075\textwidth]{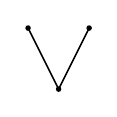}}}$ & $\;\sum b_i\alpha_{ij}\alpha_{ik}= \frac{1}{3}$\\
6. & $\vcenter{\hbox{\includegraphics[width = 0.075\textwidth]{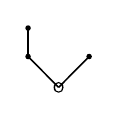}}}$ & $\;\sum b_iw_{ij}\alpha_{jk}\alpha_{jl}\beta_{lm}= \frac{1}{2}$\\
7. & $\vcenter{\hbox{\includegraphics[width = 0.075\textwidth]{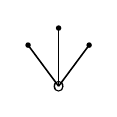}}}$ & $\;\sum b_iw_{ij}\alpha_{jk}\alpha_{jl}\alpha_{jm}= \frac{1}{1}$\\
8. & $\vcenter{\hbox{\includegraphics[width = 0.075\textwidth]{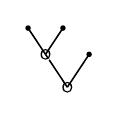}}}$ & $\;\sum b_iw_{ij}\alpha_{jk}\alpha_{jl}w_{lm}\alpha_{mn}\alpha_{mo}= \frac{1}{1}$\\
9. & $\vcenter{\hbox{\includegraphics[width = 0.075\textwidth]{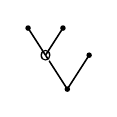}}}$ & $\;\sum b_i\alpha_{ij}\alpha_{ik}w_{kl}\alpha_{lm}\alpha_{ln}= \frac{1}{4}$\\
10. & $\vcenter{\hbox{\includegraphics[width = 0.075\textwidth]{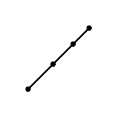}}}$ & $\;\sum b_i\beta_{ij}\beta_{jk}\beta_{kl}= \frac{1}{24}$\\
11. & $\vcenter{\hbox{\includegraphics[width = 0.075\textwidth]{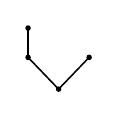}}}$ & $\;\sum b_i\alpha_{ij}\alpha_{ik}\beta_{kl}= \frac{1}{8}$\\
12. & $\vcenter{\hbox{\includegraphics[width = 0.075\textwidth]{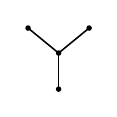}}}$ & $\;\sum b_i\beta_{ij}\alpha_{jk}\alpha_{jl}= \frac{1}{12}$\\
13. & $\vcenter{\hbox{\includegraphics[width = 0.075\textwidth]{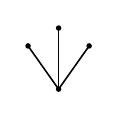}}}$ & $\;\sum b_i\alpha_{ij}\alpha_{ik}\alpha_{il}= \frac{1}{4}$\\
14. & $\vcenter{\hbox{\includegraphics[width = 0.075\textwidth]{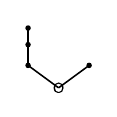}}}$ & $\;\sum b_iw_{ij}\alpha_{jk}\alpha_{jl}\beta_{lm}\beta_{mn}= \frac{1}{6}$\\
15. & $\vcenter{\hbox{\includegraphics[width = 0.075\textwidth]{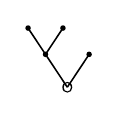}}}$ & $\;\sum b_iw_{ij}\alpha_{jk}\alpha_{jl}\alpha_{lm}\alpha_{ln}= \frac{1}{3}$\\
16. & $\vcenter{\hbox{\includegraphics[width = 0.075\textwidth]{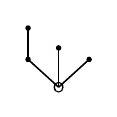}}}$ & $\;\sum b_iw_{ij}\alpha_{jk}\alpha_{jl}\alpha_{jm}\beta_{mn}= \frac{1}{2}$\\
17. & $\vcenter{\hbox{\includegraphics[width = 0.075\textwidth]{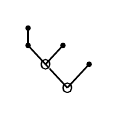}}}$ & $\;\sum b_iw_{ij}\alpha_{jk}\alpha_{jl}w_{lm}\alpha_{mn}\alpha_{mo}\beta_{op}= \frac{1}{2}$\\
18. & $\vcenter{\hbox{\includegraphics[width = 0.075\textwidth]{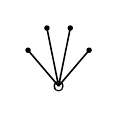}}}$ & $\;\sum b_iw_{ij}\alpha_{jk}\alpha_{jl}\alpha_{jm}\alpha_{jn}= \frac{1}{1}$\\
19. & $\vcenter{\hbox{\includegraphics[width = 0.075\textwidth]{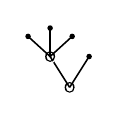}}}$ & $\;\sum b_iw_{ij}\alpha_{jk}\alpha_{jl}w_{lm}\alpha_{mn}\alpha_{mo}\alpha_{mp}= \frac{1}{1}$\\
\hline
\end{tabular}
\newpage
\begin{tabular}{lcl}
\hline
20. & $\vcenter{\hbox{\includegraphics[width = 0.075\textwidth]{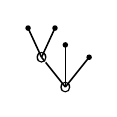}}}$ & $\;\sum b_iw_{ij}\alpha_{jk}\alpha_{jl}\alpha_{jm}w_{mn}\alpha_{no}\alpha_{np}= \frac{1}{1}$\\
21. & $\vcenter{\hbox{\includegraphics[width = 0.075\textwidth]{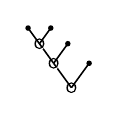}}}$ & $\;\sum b_iw_{ij}\alpha_{jk}\alpha_{jl}w_{lm}\alpha_{mn}\alpha_{mo}w_{op}\alpha_{pq}\alpha_{pr}= \frac{1}{1}$\\
22. & $\vcenter{\hbox{\includegraphics[width = 0.075\textwidth]{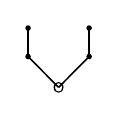}}}$ & $\;\sum b_iw_{ij}\alpha_{jk}\beta_{kl}\alpha_{jl}\beta_{lm}= \frac{1}{4}$\\
23. & $\vcenter{\hbox{\includegraphics[width = 0.075\textwidth]{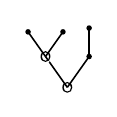}}}$ & $\;\sum b_iw_{ij}\alpha_{jk}\beta_{kl}\alpha_{jl}w_{lm}\alpha_{mn}\alpha_{mo}= \frac{1}{2}$\\
24. & $\vcenter{\hbox{\includegraphics[width = 0.075\textwidth]{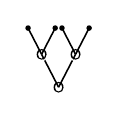}}}$ & $\;\sum b_iw_{ij}\alpha_{jk}w_{kl}\alpha_{lm}\alpha_{ln}\alpha_{jl}w_{lm}\alpha_{mn}\alpha_{mo}= \frac{1}{1}$\\
25. & $\vcenter{\hbox{\includegraphics[width = 0.075\textwidth]{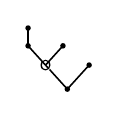}}}$ & $\;\sum b_i\alpha_{ij}\alpha_{ik}w_{kl}\alpha_{lm}\alpha_{ln}\beta_{no}= \frac{1}{10}$\\
26. & $\vcenter{\hbox{\includegraphics[width = 0.075\textwidth]{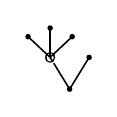}}}$ & $\;\sum b_i\alpha_{ij}\alpha_{ik}w_{kl}\alpha_{lm}\alpha_{ln}\alpha_{lo}= \frac{1}{5}$\\
27. & $\vcenter{\hbox{\includegraphics[width = 0.075\textwidth]{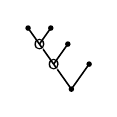}}}$ & $\;\sum b_i\alpha_{ij}\alpha_{ik}w_{kl}\alpha_{lm}\alpha_{ln}w_{no}\alpha_{op}\alpha_{oq}= \frac{1}{5}$\\
28. & $\vcenter{\hbox{\includegraphics[width = 0.075\textwidth]{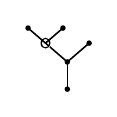}}}$ & $\;\sum b_i\beta_{ij}\alpha_{jk}\alpha_{jl}w_{lm}\alpha_{mn}\alpha_{mo}= \frac{1}{20}$\\
29. & $\vcenter{\hbox{\includegraphics[width = 0.075\textwidth]{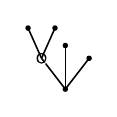}}}$ & $\;\sum b_i\alpha_{ij}\alpha_{ik}\alpha_{il}w_{lm}\alpha_{mn}\alpha_{mo}= \frac{1}{5}$\\
30. & $\vcenter{\hbox{\includegraphics[width = 0.075\textwidth]{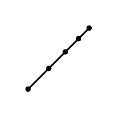}}}$ & $\;\sum b_i\beta_{ij}\beta_{jk}\beta_{kl}\beta_{lm}= \frac{1}{120}$\\
31. & $\vcenter{\hbox{\includegraphics[width = 0.075\textwidth]{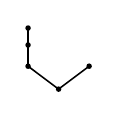}}}$ & $\;\sum b_i\alpha_{ij}\alpha_{ik}\beta_{kl}\beta_{lm}= \frac{1}{30}$\\
32. & $\vcenter{\hbox{\includegraphics[width = 0.075\textwidth]{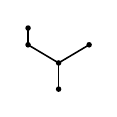}}}$ & $\;\sum b_i\beta_{ij}\alpha_{jk}\alpha_{jl}\beta_{lm}= \frac{1}{40}$\\
33. & $\vcenter{\hbox{\includegraphics[width = 0.075\textwidth]{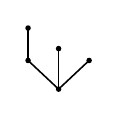}}}$ & $\;\sum b_i\alpha_{ij}\alpha_{ik}\alpha_{il}\beta_{lm}= \frac{1}{10}$\\
34. & $\vcenter{\hbox{\includegraphics[width = 0.075\textwidth]{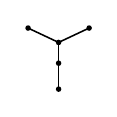}}}$ & $\;\sum b_i\beta_{ij}\beta_{jk}\alpha_{kl}\alpha_{km}= \frac{1}{60}$\\
35. & $\vcenter{\hbox{\includegraphics[width = 0.075\textwidth]{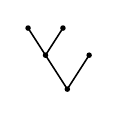}}}$ & $\;\sum b_i\alpha_{ij}\alpha_{ik}\alpha_{kl}\alpha_{km}= \frac{1}{15}$\\
36. & $\vcenter{\hbox{\includegraphics[width = 0.075\textwidth]{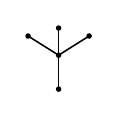}}}$ & $\;\sum b_i\beta_{ij}\alpha_{jk}\alpha_{jl}\alpha_{jm}= \frac{1}{20}$\\
37. & $\vcenter{\hbox{\includegraphics[width = 0.075\textwidth]{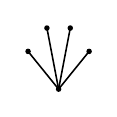}}}$ & $\;\sum b_i\alpha_{ij}\alpha_{ik}\alpha_{il}\alpha_{im}= \frac{1}{5}$\\
38. & $\vcenter{\hbox{\includegraphics[width = 0.075\textwidth]{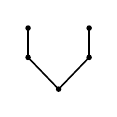}}}$ & $\;\sum b_i\alpha_{ij}\beta_{jk}\alpha_{ik}\beta_{kl}= \frac{1}{20}$\\
39. & $\vcenter{\hbox{\includegraphics[width = 0.075\textwidth]{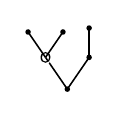}}}$ & $\;\sum b_i\alpha_{ij}\beta_{jk}\alpha_{ik}w_{kl}\alpha_{lm}\alpha_{ln}= \frac{1}{10}$\\
40. & $\vcenter{\hbox{\includegraphics[width = 0.075\textwidth]{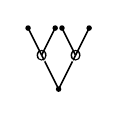}}}$ & $\;\sum b_i\alpha_{ij}w_{jk}\alpha_{kl}\alpha_{km}\alpha_{ik}w_{kl}\alpha_{lm}\alpha_{ln}= \frac{1}{5}$\\
41. & $\vcenter{\hbox{\includegraphics[width = 0.075\textwidth]{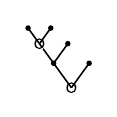}}}$ & $\;\sum b_iw_{ij}\alpha_{jk}\alpha_{jl}\alpha_{lm}\alpha_{ln}w_{no}\alpha_{op}\alpha_{oq}= \frac{1}{4}$\\
42. & $\vcenter{\hbox{\includegraphics[width = 0.075\textwidth]{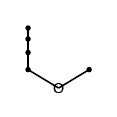}}}$ & $\;\sum b_iw_{ij}\alpha_{jk}\alpha_{jl}\beta_{lm}\beta_{mn}\beta_{no}= \frac{1}{24}$\\
43. & $\vcenter{\hbox{\includegraphics[width = 0.075\textwidth]{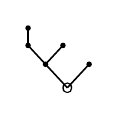}}}$ & $\;\sum b_iw_{ij}\alpha_{jk}\alpha_{jl}\alpha_{lm}\alpha_{ln}\beta_{no}= \frac{1}{8}$\\
44. & $\vcenter{\hbox{\includegraphics[width = 0.075\textwidth]{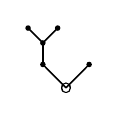}}}$ & $\;\sum b_iw_{ij}\alpha_{jk}\alpha_{jl}\beta_{lm}\alpha_{mn}\alpha_{mo}= \frac{1}{12}$\\
45. & $\vcenter{\hbox{\includegraphics[width = 0.075\textwidth]{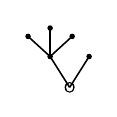}}}$ & $\;\sum b_iw_{ij}\alpha_{jk}\alpha_{jl}\alpha_{lm}\alpha_{ln}\alpha_{lo}= \frac{1}{4}$\\
\hline
\end{tabular}
\newpage
\begin{tabular}{lcl}
\hline
46. & $\vcenter{\hbox{\includegraphics[width = 0.075\textwidth]{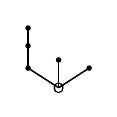}}}$ & $\;\sum b_iw_{ij}\alpha_{jk}\alpha_{jl}\alpha_{jm}\beta_{mn}\beta_{no}= \frac{1}{6}$\\
47. & $\vcenter{\hbox{\includegraphics[width = 0.075\textwidth]{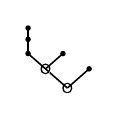}}}$ & $\;\sum b_iw_{ij}\alpha_{jk}\alpha_{jl}w_{lm}\alpha_{mn}\alpha_{mo}\beta_{op}\beta_{pq}= \frac{1}{6}$\\
48. & $\vcenter{\hbox{\includegraphics[width = 0.075\textwidth]{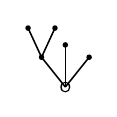}}}$ & $\;\sum b_iw_{ij}\alpha_{jk}\alpha_{jl}\alpha_{jm}\alpha_{mn}\alpha_{mo}= \frac{1}{3}$\\
49. & $\vcenter{\hbox{\includegraphics[width = 0.075\textwidth]{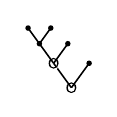}}}$ & $\;\sum b_iw_{ij}\alpha_{jk}\alpha_{jl}w_{lm}\alpha_{mn}\alpha_{mo}\alpha_{op}\alpha_{oq}= \frac{1}{3}$\\
50. & $\vcenter{\hbox{\includegraphics[width = 0.075\textwidth]{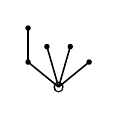}}}$ & $\;\sum b_iw_{ij}\alpha_{jk}\alpha_{jl}\alpha_{jm}\alpha_{jn}\beta_{no}= \frac{1}{2}$\\
51. & $\vcenter{\hbox{\includegraphics[width = 0.075\textwidth]{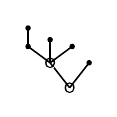}}}$ & $\;\sum b_iw_{ij}\alpha_{jk}\alpha_{jl}w_{lm}\alpha_{mn}\alpha_{mo}\alpha_{mp}\beta_{pq}= \frac{1}{2}$\\
52. & $\vcenter{\hbox{\includegraphics[width = 0.075\textwidth]{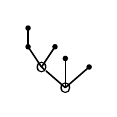}}}$ & $\;\sum b_iw_{ij}\alpha_{jk}\alpha_{jl}\alpha_{jm}w_{mn}\alpha_{no}\alpha_{np}\beta_{pq}= \frac{1}{2}$\\
53. & $\vcenter{\hbox{\includegraphics[width = 0.075\textwidth]{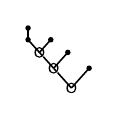}}}$ & $\;\sum b_iw_{ij}\alpha_{jk}\alpha_{jl}w_{lm}\alpha_{mn}\alpha_{mo}w_{op}\alpha_{pq}\alpha_{pr}\beta_{rs}= \frac{1}{2}$\\
54. & $\vcenter{\hbox{\includegraphics[width = 0.075\textwidth]{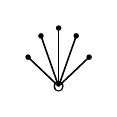}}}$ & $\;\sum b_iw_{ij}\alpha_{jk}\alpha_{jl}\alpha_{jm}\alpha_{jn}\alpha_{jo}= \frac{1}{1}$\\
55. & $\vcenter{\hbox{\includegraphics[width = 0.075\textwidth]{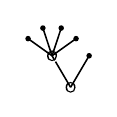}}}$ & $\;\sum b_iw_{ij}\alpha_{jk}\alpha_{jl}w_{lm}\alpha_{mn}\alpha_{mo}\alpha_{mp}\alpha_{mq}= \frac{1}{1}$\\
56. & $\vcenter{\hbox{\includegraphics[width = 0.075\textwidth]{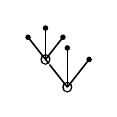}}}$ & $\;\sum b_iw_{ij}\alpha_{jk}\alpha_{jl}\alpha_{jm}w_{mn}\alpha_{no}\alpha_{np}\alpha_{nq}= \frac{1}{1}$\\
57. & $\vcenter{\hbox{\includegraphics[width = 0.075\textwidth]{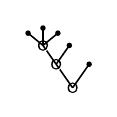}}}$ & $\;\sum b_iw_{ij}\alpha_{jk}\alpha_{jl}w_{lm}\alpha_{mn}\alpha_{mo}w_{op}\alpha_{pq}\alpha_{pr}\alpha_{ps}= \frac{1}{1}$\\
58. & $\vcenter{\hbox{\includegraphics[width = 0.075\textwidth]{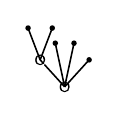}}}$ & $\;\sum b_iw_{ij}\alpha_{jk}\alpha_{jl}\alpha_{jm}\alpha_{jn}w_{no}\alpha_{op}\alpha_{oq}= \frac{1}{1}$\\
59. & $\vcenter{\hbox{\includegraphics[width = 0.075\textwidth]{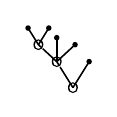}}}$ & $\;\sum b_iw_{ij}\alpha_{jk}\alpha_{jl}w_{lm}\alpha_{mn}\alpha_{mo}\alpha_{mp}w_{pq}\alpha_{qr}\alpha_{qs}= \frac{1}{1}$\\
60. & $\vcenter{\hbox{\includegraphics[width = 0.075\textwidth]{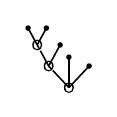}}}$ & $\;\sum b_iw_{ij}\alpha_{jk}\alpha_{jl}\alpha_{jm}w_{mn}\alpha_{no}\alpha_{np}w_{pq}\alpha_{qr}\alpha_{qs}= \frac{1}{1}$\\
61. & $\vcenter{\hbox{\includegraphics[width = 0.075\textwidth]{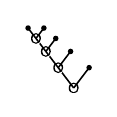}}}$ & $\;\sum b_iw_{ij}\alpha_{jk}\alpha_{jl}w_{lm}\alpha_{mn}\alpha_{mo}w_{op}\alpha_{pq}\alpha_{pr}w_{rs}\alpha_{st}\alpha_{su}= \frac{1}{1}$\\
62. & $\vcenter{\hbox{\includegraphics[width = 0.075\textwidth]{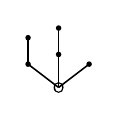}}}$ & $\;\sum b_iw_{ij}\alpha_{jk}\alpha_{jl}\beta_{lm}\alpha_{jm}\beta_{mn}= \frac{1}{4}$\\
63. & $\vcenter{\hbox{\includegraphics[width = 0.075\textwidth]{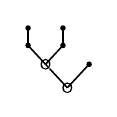}}}$ & $\;\sum b_iw_{ij}\alpha_{jk}\alpha_{jl}w_{lm}\alpha_{mn}\beta_{no}\alpha_{mo}\beta_{op}= \frac{1}{4}$\\
64. & $\vcenter{\hbox{\includegraphics[width = 0.075\textwidth]{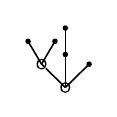}}}$ & $\;\sum b_iw_{ij}\alpha_{jk}\alpha_{jl}\beta_{lm}\alpha_{jm}w_{mn}\alpha_{no}\alpha_{np}= \frac{1}{2}$\\
65. & $\vcenter{\hbox{\includegraphics[width = 0.075\textwidth]{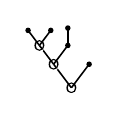}}}$ & $\;\sum b_iw_{ij}\alpha_{jk}\alpha_{jl}w_{lm}\alpha_{mn}\beta_{no}\alpha_{mo}w_{op}\alpha_{pq}\alpha_{pr}= \frac{1}{2}$\\
66. & $\vcenter{\hbox{\includegraphics[width = 0.075\textwidth]{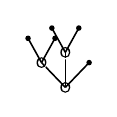}}}$ & $\;\sum b_iw_{ij}\alpha_{jk}\alpha_{jl}w_{lm}\alpha_{mn}\alpha_{mo}\alpha_{jm}w_{mn}\alpha_{no}\alpha_{np}= \frac{1}{1}$\\
67. & $\vcenter{\hbox{\includegraphics[width = 0.075\textwidth]{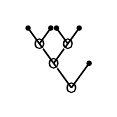}}}$ & $\;\sum b_iw_{ij}\alpha_{jk}\alpha_{jl}w_{lm}\alpha_{mn}w_{no}\alpha_{op}\alpha_{oq}\alpha_{mo}w_{op}\alpha_{pq}\alpha_{pr}= \frac{1}{1}$\\
68. & $\vcenter{\hbox{\includegraphics[width = 0.075\textwidth]{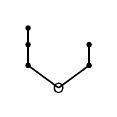}}}$ & $\;\sum b_iw_{ij}\alpha_{jk}\beta_{kl}\alpha_{jl}\beta_{lm}\beta_{mn}= \frac{1}{12}$\\
69. & $\vcenter{\hbox{\includegraphics[width = 0.075\textwidth]{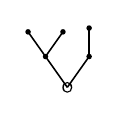}}}$ & $\;\sum b_iw_{ij}\alpha_{jk}\beta_{kl}\alpha_{jl}\alpha_{lm}\alpha_{ln}= \frac{1}{6}$\\
70. & $\vcenter{\hbox{\includegraphics[width = 0.075\textwidth]{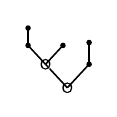}}}$ & $\;\sum b_iw_{ij}\alpha_{jk}\beta_{kl}\alpha_{jl}w_{lm}\alpha_{mn}\alpha_{mo}\beta_{op}= \frac{1}{4}$\\
\hline
\end{tabular}
\newpage
\begin{tabular}{lcl}
\hline
71. & $\vcenter{\hbox{\includegraphics[width = 0.075\textwidth]{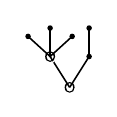}}}$ & $\;\sum b_iw_{ij}\alpha_{jk}\beta_{kl}\alpha_{jl}w_{lm}\alpha_{mn}\alpha_{mo}\alpha_{mp}= \frac{1}{2}$\\
72. & $\vcenter{\hbox{\includegraphics[width = 0.075\textwidth]{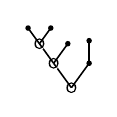}}}$ & $\;\sum b_iw_{ij}\alpha_{jk}\beta_{kl}\alpha_{jl}w_{lm}\alpha_{mn}\alpha_{mo}w_{op}\alpha_{pq}\alpha_{pr}= \frac{1}{2}$\\
73. & $\vcenter{\hbox{\includegraphics[width = 0.075\textwidth]{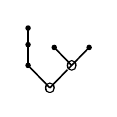}}}$ & $\;\sum b_iw_{ij}\alpha_{jk}w_{kl}\alpha_{lm}\alpha_{ln}\alpha_{jl}\beta_{lm}\beta_{mn}= \frac{1}{6}$\\
74. & $\vcenter{\hbox{\includegraphics[width = 0.075\textwidth]{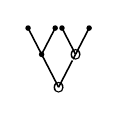}}}$ & $\;\sum b_iw_{ij}\alpha_{jk}w_{kl}\alpha_{lm}\alpha_{ln}\alpha_{jl}\alpha_{lm}\alpha_{ln}= \frac{1}{3}$\\
75. & $\vcenter{\hbox{\includegraphics[width = 0.075\textwidth]{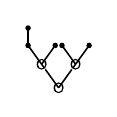}}}$ & $\;\sum b_iw_{ij}\alpha_{jk}w_{kl}\alpha_{lm}\alpha_{ln}\alpha_{jl}w_{lm}\alpha_{mn}\alpha_{mo}\beta_{op}= \frac{1}{2}$\\
76. & $\vcenter{\hbox{\includegraphics[width = 0.075\textwidth]{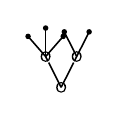}}}$ & $\;\sum b_iw_{ij}\alpha_{jk}w_{kl}\alpha_{lm}\alpha_{ln}\alpha_{jl}w_{lm}\alpha_{mn}\alpha_{mo}\alpha_{mp}= \frac{1}{1}$\\
77. & $\vcenter{\hbox{\includegraphics[width = 0.075\textwidth]{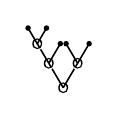}}}$ & $\;\sum b_iw_{ij}\alpha_{jk}w_{kl}\alpha_{lm}\alpha_{ln}\alpha_{jl}w_{lm}\alpha_{mn}\alpha_{mo}w_{op}\alpha_{pq}\alpha_{pr}= \frac{1}{1}$\\
78. & $\vcenter{\hbox{\includegraphics[width = 0.075\textwidth]{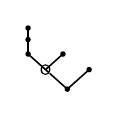}}}$ & $\;\sum b_i\alpha_{ij}\alpha_{ik}w_{kl}\alpha_{lm}\alpha_{ln}\beta_{no}\beta_{op}= \frac{1}{36}$\\
79. & $\vcenter{\hbox{\includegraphics[width = 0.075\textwidth]{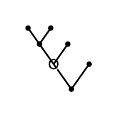}}}$ & $\;\sum b_i\alpha_{ij}\alpha_{ik}w_{kl}\alpha_{lm}\alpha_{ln}\alpha_{no}\alpha_{np}= \frac{1}{18}$\\80. & $\vcenter{\hbox{\includegraphics[width = 0.075\textwidth]{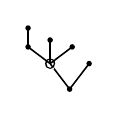}}}$ & $\;\sum b_i\alpha_{ij}\alpha_{ik}w_{kl}\alpha_{lm}\alpha_{ln}\alpha_{lo}\beta_{op}= \frac{1}{12}$\\
81. & $\vcenter{\hbox{\includegraphics[width = 0.075\textwidth]{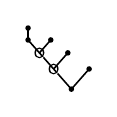}}}$ & $\;\sum b_i\alpha_{ij}\alpha_{ik}w_{kl}\alpha_{lm}\alpha_{ln}w_{no}\alpha_{op}\alpha_{oq}\beta_{qr}= \frac{1}{12}$\\
82. & $\vcenter{\hbox{\includegraphics[width = 0.075\textwidth]{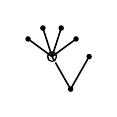}}}$ & $\;\sum b_i\alpha_{ij}\alpha_{ik}w_{kl}\alpha_{lm}\alpha_{ln}\alpha_{lo}\alpha_{lp}= \frac{1}{6}$\\
83. & $\vcenter{\hbox{\includegraphics[width = 0.075\textwidth]{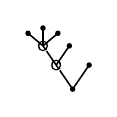}}}$ & $\;\sum b_i\alpha_{ij}\alpha_{ik}w_{kl}\alpha_{lm}\alpha_{ln}w_{no}\alpha_{op}\alpha_{oq}\alpha_{or}= \frac{1}{6}$\\
84. & $\vcenter{\hbox{\includegraphics[width = 0.075\textwidth]{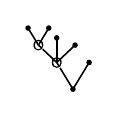}}}$ & $\;\sum b_i\alpha_{ij}\alpha_{ik}w_{kl}\alpha_{lm}\alpha_{ln}\alpha_{lo}w_{op}\alpha_{pq}\alpha_{pr}= \frac{1}{6}$\\
85. & $\vcenter{\hbox{\includegraphics[width = 0.075\textwidth]{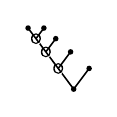}}}$ & $\;\sum b_i\alpha_{ij}\alpha_{ik}w_{kl}\alpha_{lm}\alpha_{ln}w_{no}\alpha_{op}\alpha_{oq}w_{qr}\alpha_{rs}\alpha_{rt}= \frac{1}{6}$\\
86. & $\vcenter{\hbox{\includegraphics[width = 0.075\textwidth]{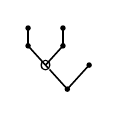}}}$ & $\;\sum b_i\alpha_{ij}\alpha_{ik}w_{kl}\alpha_{lm}\beta_{mn}\alpha_{ln}\beta_{no}= \frac{1}{24}$\\
87. & $\vcenter{\hbox{\includegraphics[width = 0.075\textwidth]{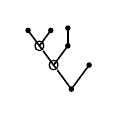}}}$ & $\;\sum b_i\alpha_{ij}\alpha_{ik}w_{kl}\alpha_{lm}\beta_{mn}\alpha_{ln}w_{no}\alpha_{op}\alpha_{oq}= \frac{1}{12}$\\
88. & $\vcenter{\hbox{\includegraphics[width = 0.075\textwidth]{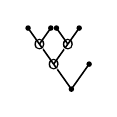}}}$ & $\;\sum b_i\alpha_{ij}\alpha_{ik}w_{kl}\alpha_{lm}w_{mn}\alpha_{no}\alpha_{np}\alpha_{ln}w_{no}\alpha_{op}\alpha_{oq}= \frac{1}{6}$\\
89. & $\vcenter{\hbox{\includegraphics[width = 0.075\textwidth]{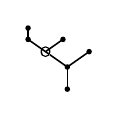}}}$ & $\;\sum b_i\beta_{ij}\alpha_{jk}\alpha_{jl}w_{lm}\alpha_{mn}\alpha_{mo}\beta_{op}= \frac{1}{60}$\\
90. & $\vcenter{\hbox{\includegraphics[width = 0.075\textwidth]{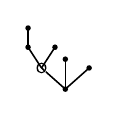}}}$ & $\;\sum b_i\alpha_{ij}\alpha_{ik}\alpha_{il}w_{lm}\alpha_{mn}\alpha_{mo}\beta_{op}= \frac{1}{12}$\\
91. & $\vcenter{\hbox{\includegraphics[width = 0.075\textwidth]{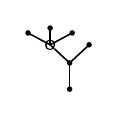}}}$ & $\;\sum b_i\beta_{ij}\alpha_{jk}\alpha_{jl}w_{lm}\alpha_{mn}\alpha_{mo}\alpha_{mp}= \frac{1}{30}$\\
92. & $\vcenter{\hbox{\includegraphics[width = 0.075\textwidth]{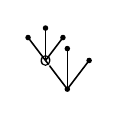}}}$ & $\;\sum b_i\alpha_{ij}\alpha_{ik}\alpha_{il}w_{lm}\alpha_{mn}\alpha_{mo}\alpha_{mp}= \frac{1}{6}$\\
93. & $\vcenter{\hbox{\includegraphics[width = 0.075\textwidth]{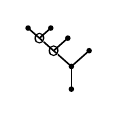}}}$ & $\;\sum b_i\beta_{ij}\alpha_{jk}\alpha_{jl}w_{lm}\alpha_{mn}\alpha_{mo}w_{op}\alpha_{pq}\alpha_{pr}= \frac{1}{30}$\\
94. & $\vcenter{\hbox{\includegraphics[width = 0.075\textwidth]{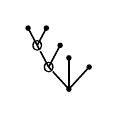}}}$ & $\;\sum b_i\alpha_{ij}\alpha_{ik}\alpha_{il}w_{lm}\alpha_{mn}\alpha_{mo}w_{op}\alpha_{pq}\alpha_{pr}= \frac{1}{6}$\\
95. & $\vcenter{\hbox{\includegraphics[width = 0.075\textwidth]{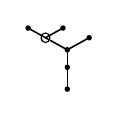}}}$ & $\;\sum b_i\beta_{ij}\beta_{jk}\alpha_{kl}\alpha_{km}w_{mn}\alpha_{no}\alpha_{np}= \frac{1}{120}$\\
\hline
\end{tabular}
\newpage
\begin{tabular}{lcl}
\hline
96. & $\vcenter{\hbox{\includegraphics[width = 0.075\textwidth]{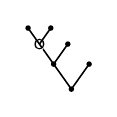}}}$ & $\;\sum b_i\alpha_{ij}\alpha_{ik}\alpha_{kl}\alpha_{km}w_{mn}\alpha_{no}\alpha_{np}= \frac{1}{24}$\\97. & $\vcenter{\hbox{\includegraphics[width = 0.075\textwidth]{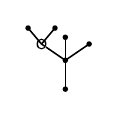}}}$ & $\;\sum b_i\beta_{ij}\alpha_{jk}\alpha_{jl}\alpha_{jm}w_{mn}\alpha_{no}\alpha_{np}= \frac{1}{30}$\\
98. & $\vcenter{\hbox{\includegraphics[width = 0.075\textwidth]{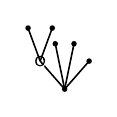}}}$ & $\;\sum b_i\alpha_{ij}\alpha_{ik}\alpha_{il}\alpha_{im}w_{mn}\alpha_{no}\alpha_{np}= \frac{1}{6}$\\
99. & $\vcenter{\hbox{\includegraphics[width = 0.075\textwidth]{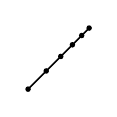}}}$ & $\;\sum b_i\beta_{ij}\beta_{jk}\beta_{kl}\beta_{lm}\beta_{mn}= \frac{1}{720}$\\
100. & $\vcenter{\hbox{\includegraphics[width = 0.075\textwidth]{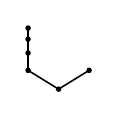}}}$ & $\;\sum b_i\alpha_{ij}\alpha_{ik}\beta_{kl}\beta_{lm}\beta_{mn}= \frac{1}{144}$\\
101. & $\vcenter{\hbox{\includegraphics[width = 0.075\textwidth]{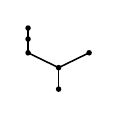}}}$ & $\;\sum b_i\beta_{ij}\alpha_{jk}\alpha_{jl}\beta_{lm}\beta_{mn}= \frac{1}{180}$\\
102. & $\vcenter{\hbox{\includegraphics[width = 0.075\textwidth]{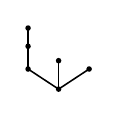}}}$ & $\;\sum b_i\alpha_{ij}\alpha_{ik}\alpha_{il}\beta_{lm}\beta_{mn}= \frac{1}{36}$\\
103. & $\vcenter{\hbox{\includegraphics[width = 0.075\textwidth]{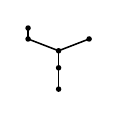}}}$ & $\;\sum b_i\beta_{ij}\beta_{jk}\alpha_{kl}\alpha_{km}\beta_{mn}= \frac{1}{240}$\\
104. & $\vcenter{\hbox{\includegraphics[width = 0.075\textwidth]{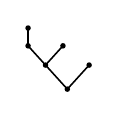}}}$ & $\;\sum b_i\alpha_{ij}\alpha_{ik}\alpha_{kl}\alpha_{km}\beta_{mn}= \frac{1}{48}$\\
105. & $\vcenter{\hbox{\includegraphics[width = 0.075\textwidth]{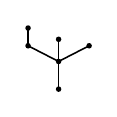}}}$ & $\;\sum b_i\beta_{ij}\alpha_{jk}\alpha_{jl}\alpha_{jm}\beta_{mn}= \frac{1}{60}$\\
106. & $\vcenter{\hbox{\includegraphics[width = 0.075\textwidth]{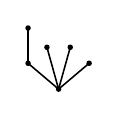}}}$ & $\;\sum b_i\alpha_{ij}\alpha_{ik}\alpha_{il}\alpha_{im}\beta_{mn}= \frac{1}{12}$\\
107. & $\vcenter{\hbox{\includegraphics[width = 0.075\textwidth]{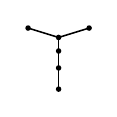}}}$ & $\;\sum b_i\beta_{ij}\beta_{jk}\beta_{kl}\alpha_{lm}\alpha_{ln}= \frac{1}{360}$\\
108. & $\vcenter{\hbox{\includegraphics[width = 0.075\textwidth]{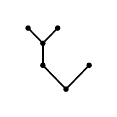}}}$ & $\;\sum b_i\alpha_{ij}\alpha_{ik}\beta_{kl}\alpha_{lm}\alpha_{ln}= \frac{1}{72}$\\
109. & $\vcenter{\hbox{\includegraphics[width = 0.075\textwidth]{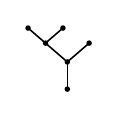}}}$ & $\;\sum b_i\beta_{ij}\alpha_{jk}\alpha_{jl}\alpha_{lm}\alpha_{ln}= \frac{1}{90}$\\
110. & $\vcenter{\hbox{\includegraphics[width = 0.075\textwidth]{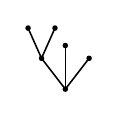}}}$ & $\;\sum b_i\alpha_{ij}\alpha_{ik}\alpha_{il}\alpha_{lm}\alpha_{ln}= \frac{1}{18}$\\
111. & $\vcenter{\hbox{\includegraphics[width = 0.075\textwidth]{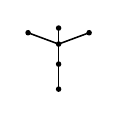}}}$ & $\;\sum b_i\beta_{ij}\beta_{jk}\alpha_{kl}\alpha_{km}\alpha_{kn}= \frac{1}{120}$\\
112. & $\vcenter{\hbox{\includegraphics[width = 0.075\textwidth]{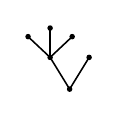}}}$ & $\;\sum b_i\alpha_{ij}\alpha_{ik}\alpha_{kl}\alpha_{km}\alpha_{kn}= \frac{1}{24}$\\
113. & $\vcenter{\hbox{\includegraphics[width = 0.075\textwidth]{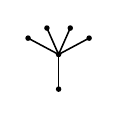}}}$ & $\;\sum b_i\beta_{ij}\alpha_{jk}\alpha_{jl}\alpha_{jm}\alpha_{jn}= \frac{1}{30}$\\
114. & $\vcenter{\hbox{\includegraphics[width = 0.075\textwidth]{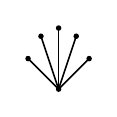}}}$ & $\;\sum b_i\alpha_{ij}\alpha_{ik}\alpha_{il}\alpha_{im}\alpha_{in}= \frac{1}{6}$\\
115. & $\vcenter{\hbox{\includegraphics[width = 0.075\textwidth]{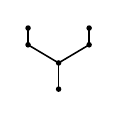}}}$ & $\;\sum b_i\beta_{ij}\alpha_{jk}\beta_{kl}\alpha_{jl}\beta_{lm}= \frac{1}{120}$\\
116. & $\vcenter{\hbox{\includegraphics[width = 0.075\textwidth]{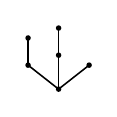}}}$ & $\;\sum b_i\alpha_{ij}\alpha_{ik}\beta_{kl}\alpha_{il}\beta_{lm}= \frac{1}{24}$\\
117. & $\vcenter{\hbox{\includegraphics[width = 0.075\textwidth]{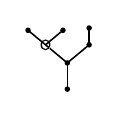}}}$ & $\;\sum b_i\beta_{ij}\alpha_{jk}\beta_{kl}\alpha_{jl}w_{lm}\alpha_{mn}\alpha_{mo}= \frac{1}{60}$\\118. & $\vcenter{\hbox{\includegraphics[width = 0.075\textwidth]{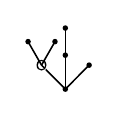}}}$ & $\;\sum b_i\alpha_{ij}\alpha_{ik}\beta_{kl}\alpha_{il}w_{lm}\alpha_{mn}\alpha_{mo}= \frac{1}{12}$\\
119. & $\vcenter{\hbox{\includegraphics[width = 0.075\textwidth]{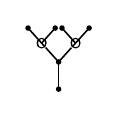}}}$ & $\;\sum b_i\beta_{ij}\alpha_{jk}w_{kl}\alpha_{lm}\alpha_{ln}\alpha_{jl}w_{lm}\alpha_{mn}\alpha_{mo}= \frac{1}{30}$\\
120. & $\vcenter{\hbox{\includegraphics[width = 0.075\textwidth]{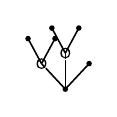}}}$ & $\;\sum b_i\alpha_{ij}\alpha_{ik}w_{kl}\alpha_{lm}\alpha_{ln}\alpha_{il}w_{lm}\alpha_{mn}\alpha_{mo}= \frac{1}{6}$\\
\hline
\end{tabular}
\newpage
\begin{tabular}{lcl}
\hline
121. & $\vcenter{\hbox{\includegraphics[width = 0.075\textwidth]{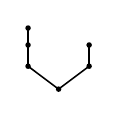}}}$ & $\;\sum b_i\alpha_{ij}\beta_{jk}\alpha_{ik}\beta_{kl}\beta_{lm}= \frac{1}{72}$\\
122. & $\vcenter{\hbox{\includegraphics[width = 0.075\textwidth]{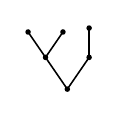}}}$ & $\;\sum b_i\alpha_{ij}\beta_{jk}\alpha_{ik}\alpha_{kl}\alpha_{km}= \frac{1}{36}$\\
123. & $\vcenter{\hbox{\includegraphics[width = 0.075\textwidth]{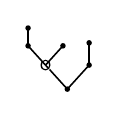}}}$ & $\;\sum b_i\alpha_{ij}\beta_{jk}\alpha_{ik}w_{kl}\alpha_{lm}\alpha_{ln}\beta_{no}= \frac{1}{24}$\\124. & $\vcenter{\hbox{\includegraphics[width = 0.075\textwidth]{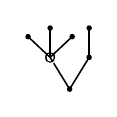}}}$ & $\;\sum b_i\alpha_{ij}\beta_{jk}\alpha_{ik}w_{kl}\alpha_{lm}\alpha_{ln}\alpha_{lo}= \frac{1}{12}$\\
125. & $\vcenter{\hbox{\includegraphics[width = 0.075\textwidth]{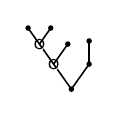}}}$ & $\;\sum b_i\alpha_{ij}\beta_{jk}\alpha_{ik}w_{kl}\alpha_{lm}\alpha_{ln}w_{no}\alpha_{op}\alpha_{oq}= \frac{1}{12}$\\
126. & $\vcenter{\hbox{\includegraphics[width = 0.075\textwidth]{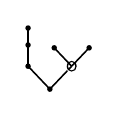}}}$ & $\;\sum b_i\alpha_{ij}w_{jk}\alpha_{kl}\alpha_{km}\alpha_{ik}\beta_{kl}\beta_{lm}= \frac{1}{36}$\\127. & $\vcenter{\hbox{\includegraphics[width = 0.075\textwidth]{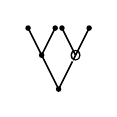}}}$ & $\;\sum b_i\alpha_{ij}w_{jk}\alpha_{kl}\alpha_{km}\alpha_{ik}\alpha_{kl}\alpha_{km}= \frac{1}{18}$\\
128. & $\vcenter{\hbox{\includegraphics[width = 0.075\textwidth]{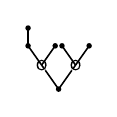}}}$ & $\;\sum b_i\alpha_{ij}w_{jk}\alpha_{kl}\alpha_{km}\alpha_{ik}w_{kl}\alpha_{lm}\alpha_{ln}\beta_{no}= \frac{1}{12}$\\
129. & $\vcenter{\hbox{\includegraphics[width = 0.075\textwidth]{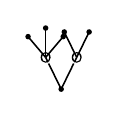}}}$ & $\;\sum b_i\alpha_{ij}w_{jk}\alpha_{kl}\alpha_{km}\alpha_{ik}w_{kl}\alpha_{lm}\alpha_{ln}\alpha_{lo}= \frac{1}{6}$\\
130. & $\vcenter{\hbox{\includegraphics[width = 0.075\textwidth]{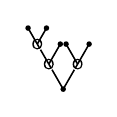}}}$ & $\;\sum b_i\alpha_{ij}w_{jk}\alpha_{kl}\alpha_{km}\alpha_{ik}w_{kl}\alpha_{lm}\alpha_{ln}w_{no}\alpha_{op}\alpha_{oq}= \frac{1}{6}$\\
\hline
\end{tabular}


\newpage
\begin{table}[H] 
\caption{Order conditions for a Rosenbrock method of type (\ref{eq_rowda1},\ref{eq_rowda2},\ref{eq_rowda3}) up to order $p=5$.}
\label{tab:tsit5da}
\end{table}
\vspace*{-0.5cm}
\begin{tabular}{lcl}
\hline
1. & $\vcenter{\hbox{\includegraphics[width = 0.075\textwidth]{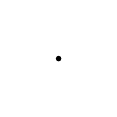}}}$ & $\;\sum b_i= \frac{1}{1}$\\
2. & $\vcenter{\hbox{\includegraphics[width = 0.075\textwidth]{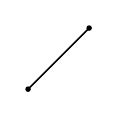}}}$ & $\;\sum b_i\alpha_{ij}= \frac{1}{2}$\\
3. & $\vcenter{\hbox{\includegraphics[width = 0.075\textwidth]{baum_3.pdf}}}$ & $\;\sum b_iw_{ij}\alpha_{jk}\alpha_{jl}= \frac{1}{1}$\\
4. & $\vcenter{\hbox{\includegraphics[width = 0.075\textwidth]{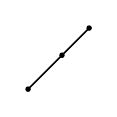}}}$ & $\;\sum b_i\alpha_{ij}\alpha_{jk}= \frac{1}{6}$\\
5. & $\vcenter{\hbox{\includegraphics[width = 0.075\textwidth]{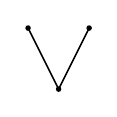}}}$ & $\;\sum b_i\alpha_{ij}\alpha_{ik}= \frac{1}{3}$\\
6. & $\vcenter{\hbox{\includegraphics[width = 0.075\textwidth]{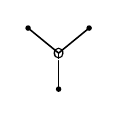}}}$ & $\;\sum b_i\alpha_{ij}w_{jk}\alpha_{kl}\alpha_{km}= \frac{1}{3}$\\
7. & $\vcenter{\hbox{\includegraphics[width = 0.075\textwidth]{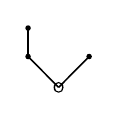}}}$ & $\;\sum b_iw_{ij}\alpha_{jk}\alpha_{jl}\alpha_{lm}= \frac{1}{2}$\\
8. & $\vcenter{\hbox{\includegraphics[width = 0.075\textwidth]{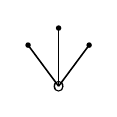}}}$ & $\;\sum b_iw_{ij}\alpha_{jk}\alpha_{jl}\alpha_{jm}= \frac{1}{1}$\\
9. & $\vcenter{\hbox{\includegraphics[width = 0.075\textwidth]{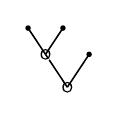}}}$ & $\;\sum b_iw_{ij}\alpha_{jk}\alpha_{jl}w_{lm}\alpha_{mn}\alpha_{mo}= \frac{1}{1}$\\
10. & $\vcenter{\hbox{\includegraphics[width = 0.075\textwidth]{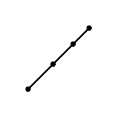}}}$ & $\;\sum b_i\alpha_{ij}\alpha_{jk}\alpha_{kl}= \frac{1}{24}$\\
11. & $\vcenter{\hbox{\includegraphics[width = 0.075\textwidth]{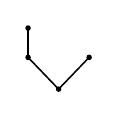}}}$ & $\;\sum b_i\alpha_{ij}\alpha_{ik}\alpha_{kl}= \frac{1}{8}$\\
12. & $\vcenter{\hbox{\includegraphics[width = 0.075\textwidth]{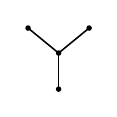}}}$ & $\;\sum b_i\alpha_{ij}\alpha_{jk}\alpha_{jl}= \frac{1}{12}$\\
13. & $\vcenter{\hbox{\includegraphics[width = 0.075\textwidth]{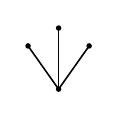}}}$ & $\;\sum b_i\alpha_{ij}\alpha_{ik}\alpha_{il}= \frac{1}{4}$\\
14. & $\vcenter{\hbox{\includegraphics[width = 0.075\textwidth]{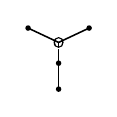}}}$ & $\;\sum b_i\alpha_{ij}\alpha_{jk}w_{kl}\alpha_{lm}\alpha_{ln}= \frac{1}{12}$\\
15. & $\vcenter{\hbox{\includegraphics[width = 0.075\textwidth]{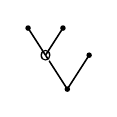}}}$ & $\;\sum b_i\alpha_{ij}\alpha_{ik}w_{kl}\alpha_{lm}\alpha_{ln}= \frac{1}{4}$\\
16. & $\vcenter{\hbox{\includegraphics[width = 0.075\textwidth]{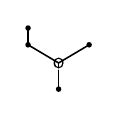}}}$ & $\;\sum b_i\alpha_{ij}w_{jk}\alpha_{kl}\alpha_{km}\alpha_{mn}= \frac{1}{8}$\\
17. & $\vcenter{\hbox{\includegraphics[width = 0.075\textwidth]{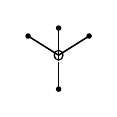}}}$ & $\;\sum b_i\alpha_{ij}w_{jk}\alpha_{kl}\alpha_{km}\alpha_{kn}= \frac{1}{4}$\\
18. & $\vcenter{\hbox{\includegraphics[width = 0.075\textwidth]{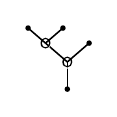}}}$ & $\;\sum b_i\alpha_{ij}w_{jk}\alpha_{kl}\alpha_{km}w_{mn}\alpha_{no}\alpha_{np}= \frac{1}{4}$\\
19. & $\vcenter{\hbox{\includegraphics[width = 0.075\textwidth]{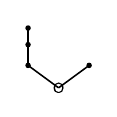}}}$ & $\;\sum b_iw_{ij}\alpha_{jk}\alpha_{jl}\alpha_{lm}\alpha_{mn}= \frac{1}{6}$\\
20. & $\vcenter{\hbox{\includegraphics[width = 0.075\textwidth]{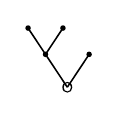}}}$ & $\;\sum b_iw_{ij}\alpha_{jk}\alpha_{jl}\alpha_{lm}\alpha_{ln}= \frac{1}{3}$\\
\hline
\end{tabular}
\newpage
\begin{tabular}{lcl}
\hline
21. & $\vcenter{\hbox{\includegraphics[width = 0.075\textwidth]{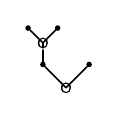}}}$ & $\;\sum b_iw_{ij}\alpha_{jk}\alpha_{jl}\alpha_{lm}w_{mn}\alpha_{no}\alpha_{np}= \frac{1}{3}$\\
22. & $\vcenter{\hbox{\includegraphics[width = 0.075\textwidth]{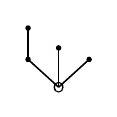}}}$ & $\;\sum b_iw_{ij}\alpha_{jk}\alpha_{jl}\alpha_{jm}\alpha_{mn}= \frac{1}{2}$\\
23. & $\vcenter{\hbox{\includegraphics[width = 0.075\textwidth]{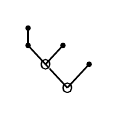}}}$ & $\;\sum b_iw_{ij}\alpha_{jk}\alpha_{jl}w_{lm}\alpha_{mn}\alpha_{mo}\alpha_{op}= \frac{1}{2}$\\
24. & $\vcenter{\hbox{\includegraphics[width = 0.075\textwidth]{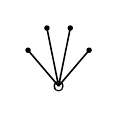}}}$ & $\;\sum b_iw_{ij}\alpha_{jk}\alpha_{jl}\alpha_{jm}\alpha_{jn}= \frac{1}{1}$\\
25. & $\vcenter{\hbox{\includegraphics[width = 0.075\textwidth]{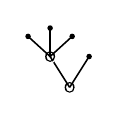}}}$ & $\;\sum b_iw_{ij}\alpha_{jk}\alpha_{jl}w_{lm}\alpha_{mn}\alpha_{mo}\alpha_{mp}= \frac{1}{1}$\\
26. & $\vcenter{\hbox{\includegraphics[width = 0.075\textwidth]{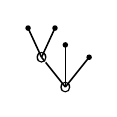}}}$ & $\;\sum b_iw_{ij}\alpha_{jk}\alpha_{jl}\alpha_{jm}w_{mn}\alpha_{no}\alpha_{np}= \frac{1}{1}$\\
27. & $\vcenter{\hbox{\includegraphics[width = 0.075\textwidth]{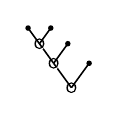}}}$ & $\;\sum b_iw_{ij}\alpha_{jk}\alpha_{jl}w_{lm}\alpha_{mn}\alpha_{mo}w_{op}\alpha_{pq}\alpha_{pr}= \frac{1}{1}$\\        
28. & $\vcenter{\hbox{\includegraphics[width = 0.075\textwidth]{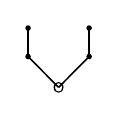}}}$ & $\;\sum b_iw_{ij}\alpha_{jk}\alpha_{kl}\alpha_{jl}\alpha_{lm}= \frac{1}{4}$\\
29. & $\vcenter{\hbox{\includegraphics[width = 0.075\textwidth]{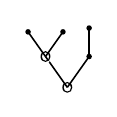}}}$ & $\;\sum b_iw_{ij}\alpha_{jk}\alpha_{kl}\alpha_{jl}w_{lm}\alpha_{mn}\alpha_{mo}= \frac{1}{2}$\\
30. & $\vcenter{\hbox{\includegraphics[width = 0.075\textwidth]{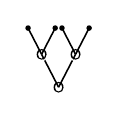}}}$ & $\;\sum b_iw_{ij}\alpha_{jk}w_{kl}\alpha_{lm}\alpha_{ln}\alpha_{jl}w_{lm}\alpha_{mn}\alpha_{mo}= \frac{1}{1}$\\        
31. & $\vcenter{\hbox{\includegraphics[width = 0.075\textwidth]{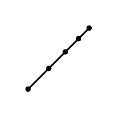}}}$ & $\;\sum b_i\alpha_{ij}\alpha_{jk}\alpha_{kl}\alpha_{lm}= \frac{1}{120}$\\
32. & $\vcenter{\hbox{\includegraphics[width = 0.075\textwidth]{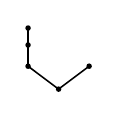}}}$ & $\;\sum b_i\alpha_{ij}\alpha_{ik}\alpha_{kl}\alpha_{lm}= \frac{1}{30}$\\
33. & $\vcenter{\hbox{\includegraphics[width = 0.075\textwidth]{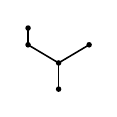}}}$ & $\;\sum b_i\alpha_{ij}\alpha_{jk}\alpha_{jl}\alpha_{lm}= \frac{1}{40}$\\
34. & $\vcenter{\hbox{\includegraphics[width = 0.075\textwidth]{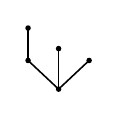}}}$ & $\;\sum b_i\alpha_{ij}\alpha_{ik}\alpha_{il}\alpha_{lm}= \frac{1}{10}$\\
35. & $\vcenter{\hbox{\includegraphics[width = 0.075\textwidth]{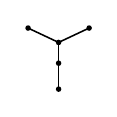}}}$ & $\;\sum b_i\alpha_{ij}\alpha_{jk}\alpha_{kl}\alpha_{km}= \frac{1}{60}$\\
36. & $\vcenter{\hbox{\includegraphics[width = 0.075\textwidth]{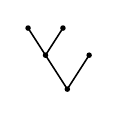}}}$ & $\;\sum b_i\alpha_{ij}\alpha_{ik}\alpha_{kl}\alpha_{km}= \frac{1}{15}$\\
37. & $\vcenter{\hbox{\includegraphics[width = 0.075\textwidth]{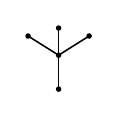}}}$ & $\;\sum b_i\alpha_{ij}\alpha_{jk}\alpha_{jl}\alpha_{jm}= \frac{1}{20}$\\
38. & $\vcenter{\hbox{\includegraphics[width = 0.075\textwidth]{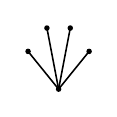}}}$ & $\;\sum b_i\alpha_{ij}\alpha_{ik}\alpha_{il}\alpha_{im}= \frac{1}{5}$\\
39. & $\vcenter{\hbox{\includegraphics[width = 0.075\textwidth]{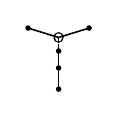}}}$ & $\;\sum b_i\alpha_{ij}\alpha_{jk}\alpha_{kl}w_{lm}\alpha_{mn}\alpha_{mo}= \frac{1}{60}$\\
40. & $\vcenter{\hbox{\includegraphics[width = 0.075\textwidth]{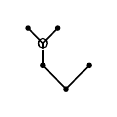}}}$ & $\;\sum b_i\alpha_{ij}\alpha_{ik}\alpha_{kl}w_{lm}\alpha_{mn}\alpha_{mo}= \frac{1}{15}$\\
41. & $\vcenter{\hbox{\includegraphics[width = 0.075\textwidth]{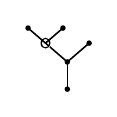}}}$ & $\;\sum b_i\alpha_{ij}\alpha_{jk}\alpha_{jl}w_{lm}\alpha_{mn}\alpha_{mo}= \frac{1}{20}$\\
42. & $\vcenter{\hbox{\includegraphics[width = 0.075\textwidth]{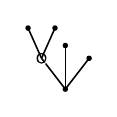}}}$ & $\;\sum b_i\alpha_{ij}\alpha_{ik}\alpha_{il}w_{lm}\alpha_{mn}\alpha_{mo}= \frac{1}{5}$\\
43. & $\vcenter{\hbox{\includegraphics[width = 0.075\textwidth]{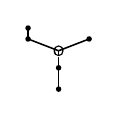}}}$ & $\;\sum b_i\alpha_{ij}\alpha_{jk}w_{kl}\alpha_{lm}\alpha_{ln}\alpha_{no}= \frac{1}{40}$\\
44. & $\vcenter{\hbox{\includegraphics[width = 0.075\textwidth]{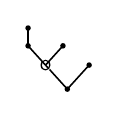}}}$ & $\;\sum b_i\alpha_{ij}\alpha_{ik}w_{kl}\alpha_{lm}\alpha_{ln}\alpha_{no}= \frac{1}{10}$\\
45. & $\vcenter{\hbox{\includegraphics[width = 0.075\textwidth]{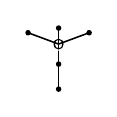}}}$ & $\;\sum b_i\alpha_{ij}\alpha_{jk}w_{kl}\alpha_{lm}\alpha_{ln}\alpha_{lo}= \frac{1}{20}$\\
\hline
\end{tabular}
\newpage
\begin{tabular}{lcl}
\hline
46. & $\vcenter{\hbox{\includegraphics[width = 0.075\textwidth]{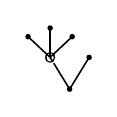}}}$ & $\;\sum b_i\alpha_{ij}\alpha_{ik}w_{kl}\alpha_{lm}\alpha_{ln}\alpha_{lo}= \frac{1}{5}$\\
47. & $\vcenter{\hbox{\includegraphics[width = 0.075\textwidth]{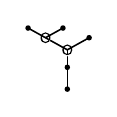}}}$ & $\;\sum b_i\alpha_{ij}\alpha_{jk}w_{kl}\alpha_{lm}\alpha_{ln}w_{no}\alpha_{op}\alpha_{oq}= \frac{1}{20}$\\
48. & $\vcenter{\hbox{\includegraphics[width = 0.075\textwidth]{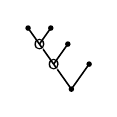}}}$ & $\;\sum b_i\alpha_{ij}\alpha_{ik}w_{kl}\alpha_{lm}\alpha_{ln}w_{no}\alpha_{op}\alpha_{oq}= \frac{1}{5}$\\
49. & $\vcenter{\hbox{\includegraphics[width = 0.075\textwidth]{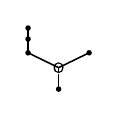}}}$ & $\;\sum b_i\alpha_{ij}w_{jk}\alpha_{kl}\alpha_{km}\alpha_{mn}\alpha_{no}= \frac{1}{30}$\\
50. & $\vcenter{\hbox{\includegraphics[width = 0.075\textwidth]{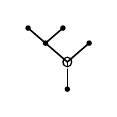}}}$ & $\;\sum b_i\alpha_{ij}w_{jk}\alpha_{kl}\alpha_{km}\alpha_{mn}\alpha_{mo}= \frac{1}{15}$\\
51. & $\vcenter{\hbox{\includegraphics[width = 0.075\textwidth]{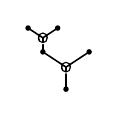}}}$ & $\;\sum b_i\alpha_{ij}w_{jk}\alpha_{kl}\alpha_{km}\alpha_{mn}w_{no}\alpha_{op}\alpha_{oq}= \frac{1}{15}$\\
52. & $\vcenter{\hbox{\includegraphics[width = 0.075\textwidth]{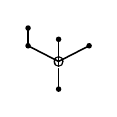}}}$ & $\;\sum b_i\alpha_{ij}w_{jk}\alpha_{kl}\alpha_{km}\alpha_{kn}\alpha_{no}= \frac{1}{10}$\\
53. & $\vcenter{\hbox{\includegraphics[width = 0.075\textwidth]{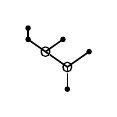}}}$ & $\;\sum b_i\alpha_{ij}w_{jk}\alpha_{kl}\alpha_{km}w_{mn}\alpha_{no}\alpha_{np}\alpha_{pq}= \frac{1}{10}$\\
54. & $\vcenter{\hbox{\includegraphics[width = 0.075\textwidth]{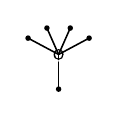}}}$ & $\;\sum b_i\alpha_{ij}w_{jk}\alpha_{kl}\alpha_{km}\alpha_{kn}\alpha_{ko}= \frac{1}{5}$\\
55. & $\vcenter{\hbox{\includegraphics[width = 0.075\textwidth]{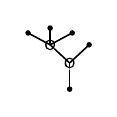}}}$ & $\;\sum b_i\alpha_{ij}w_{jk}\alpha_{kl}\alpha_{km}w_{mn}\alpha_{no}\alpha_{np}\alpha_{nq}= \frac{1}{5}$\\
56. & $\vcenter{\hbox{\includegraphics[width = 0.075\textwidth]{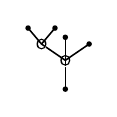}}}$ & $\;\sum b_i\alpha_{ij}w_{jk}\alpha_{kl}\alpha_{km}\alpha_{kn}w_{no}\alpha_{op}\alpha_{oq}= \frac{1}{5}$\\
57. & $\vcenter{\hbox{\includegraphics[width = 0.075\textwidth]{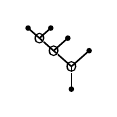}}}$ & $\;\sum b_i\alpha_{ij}w_{jk}\alpha_{kl}\alpha_{km}w_{mn}\alpha_{no}\alpha_{np}w_{pq}\alpha_{qr}\alpha_{qs}= \frac{1}{5}$\\
58. & $\vcenter{\hbox{\includegraphics[width = 0.075\textwidth]{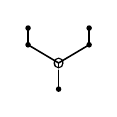}}}$ & $\;\sum b_i\alpha_{ij}w_{jk}\alpha_{kl}\alpha_{lm}\alpha_{km}\alpha_{mn}= \frac{1}{20}$\\
59. & $\vcenter{\hbox{\includegraphics[width = 0.075\textwidth]{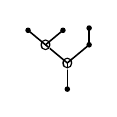}}}$ & $\;\sum b_i\alpha_{ij}w_{jk}\alpha_{kl}\alpha_{lm}\alpha_{km}w_{mn}\alpha_{no}\alpha_{np}= \frac{1}{10}$\\
60. & $\vcenter{\hbox{\includegraphics[width = 0.075\textwidth]{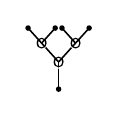}}}$ & $\;\sum b_i\alpha_{ij}w_{jk}\alpha_{kl}w_{lm}\alpha_{mn}\alpha_{mo}\alpha_{km}w_{mn}\alpha_{no}\alpha_{np}= \frac{1}{5}$\\
61. & $\vcenter{\hbox{\includegraphics[width = 0.075\textwidth]{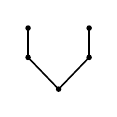}}}$ & $\;\sum b_i\alpha_{ij}\alpha_{jk}\alpha_{ik}\alpha_{kl}= \frac{1}{20}$\\
62. & $\vcenter{\hbox{\includegraphics[width = 0.075\textwidth]{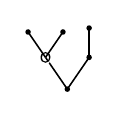}}}$ & $\;\sum b_i\alpha_{ij}\alpha_{jk}\alpha_{ik}w_{kl}\alpha_{lm}\alpha_{ln}= \frac{1}{10}$\\
63. & $\vcenter{\hbox{\includegraphics[width = 0.075\textwidth]{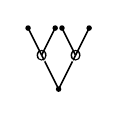}}}$ & $\;\sum b_i\alpha_{ij}w_{jk}\alpha_{kl}\alpha_{km}\alpha_{ik}w_{kl}\alpha_{lm}\alpha_{ln}= \frac{1}{5}$\\
\hline
\end{tabular}

\newpage
\begin{table}[H] 
\caption{Coefficients of method {\tt Tsit5DA}}
\label{tab:coeff}
\end{table}
\begin{verbatim}
  alpha=[0.0 0.0 0.0 0.0 0.0 0.0 0.0 0.0 0.0 0.0 0.0 0.0; 
         0.3 0.0 0.0 0.0 0.0 0.0 0.0 0.0 0.0 0.0 0.0 0.0; 
         0.4 0.0 0.0 0.0 0.0 0.0 0.0 0.0 0.0 0.0 0.0 0.0; 
         0.161 0.0 0.0 0.0 0.0 0.0 0.0 0.0 0.0 0.0 0.0 0.0; 
         -0.008480655492356989 0.0 0.0 0.335480655492357 0.0 0.0 
           0.0 0.0 0.0 0.0 0.0 0.0; 
         2.8971530571054935 0.0 0.0 -6.359448489975075 4.3622954328695815 0.0 
           0.0 0.0 0.0 0.0 0.0 0.0; 
         5.325864828439257 0.0 0.0 -11.748883564062828 7.4955393428898365 -0.09249506636175525 
           0.0 0.0 0.0 0.0 0.0 0.0; 
         5.86145544294642 0.0 0.0 -12.92096931784711 8.159367898576159 -0.071584973281401 
           -0.028269050394068383 0.0 0.0 0.0 0.0 0.0; 
         0.09646076681806523 0.0 0.0 0.01 0.4798896504144996 1.379008574103742 
           -3.290069515436081 2.324710524099774 0.0 0.0 0.0 0.0; 
         0.09468075576583945 0.0 0.0 0.009183565540343254 0.4877705284247616 1.234297566930479 
           -2.7077123499835256 1.866628418170587 0.015151515151515152 0.0 0.0 0.0; 
         0.09646076681806523 0.0 0.0 0.01 0.4798896504144996 1.379008574103742 
           -3.290069515436081 2.324710524099774 0.0 0.0 0.0 0.0; 
         0.09468075576583945 0.0 0.0 0.009183565540343254 0.4877705284247616 1.234297566930479 
           -2.7077123499835256 1.866628418170587 -0.13484848484848483 0.0 0.15 0.0]
  gamma=[0.15 0.0 0.0 0.0 0.0 0.0 0.0 0.0 0.0 0.0 0.0 0.0; 
         0.5470689774431368 0.15 0.0 0.0 0.0 0.0 0.0 0.0 0.0 0.0 0.0 0.0; 
         -0.0723537422175421 0.0666666666666667 0.15 0.0 0.0 0.0 0.0 0.0 0.0 0.0 0.0 0.0; 
         -0.11997574346406034 -0.20497635844374418 0.1257585188328081 0.15 0.0 0.0 
           0.0 0.0 0.0 0.0 0.0 0.0; 
         0.3751214208728726 -0.6896518858336065 0.355777003175544 0.09308620463102296 0.15 0.0 
           0.0 0.0 0.0 0.0 0.0 0.0; 
         -2.339423457351162 -1.8924202822866893 1.3476713525236836 7.143916166630147 -3.8352059902547007 
           0.15 0.0 0.0 0.0 0.0 0.0 0.0; 
         -4.632327787862374 -0.9275563213580595 1.3114822266754764 12.288465257549579 -7.550172308571812 
           0.11237010207373185 0.15 0.0 0.0 0.0 0.0 0.0; 
         -5.308384000531637 -1.235796359903477 1.4327893840055572 13.611173348816065 -8.203424318957262 
           0.23478742833475824 -0.06966253474809248 0.15 0.0 0.0 0.0 0.0; 
         0.6035096617978578 3.7030920005107406 9.236101686975612 1.1223090015867678 -8.707588403514192 
           -10.01583191268519 3.226138565592647 3.563871912389068 0.15 0.0 0.0 0.0; 
         0.5358920454864625 0.5149989566328188 -2.906166595272873 0.28758667283221606 0.4409793917839428 
           -1.2462207699816854 2.8597299754852776 -1.7759657086671305 0.7624212212647992 0.15 0.0 0.0; 
         -0.0017800110522257773 0.0 0.0 -0.0008164344596567463 0.007880878010261994 -0.1447110071732629 
           0.5823571654525552 -0.45808210592918686 -0.13484848484848483 0.0 0.15 0.0; 
         0.0017800110522257773 0.0 0.0 0.0008164344596567463 -0.007880878010261994 0.1447110071732629 
           -0.5823571654525552 0.45808210592918686 0.13484848484848483 -0.15 -0.15 0.15]
  b=[0.09646076681806523, 0.0, 0.0, 0.01, 0.4798896504144996, 1.379008574103742, 
     -3.290069515436081, 2.324710524099774, 0.0, -0.15, 0.0, 0.15]
  bd=[0.09468075576583945, 0.0, 0.0, 0.009183565540343254, 0.4877705284247616, 1.234297566930479, 
      -2.7077123499835256, 1.866628418170587, -0.13484848484848483, 0.0, 0.15, 0.0]
  c=[-0.8556749116393667, 0.1165263061110306, -0.038120922841221455, -0.15789728749504028, 
     0.54499490500098, 1.0853086321284309, -2.2958098031370873, 1.566895939698076,  
     8.34587614295097, -0.4162190065087707, -8.314552638841711, 0.41867264457370923]
  d=[5.79723517059224, 9.361429135834928, -3.062538663421373, -13.568052287784441, 
     1.3736819148585004, -2.344366172070166, 9.053170825304539, -7.042985092806263, 
     -147.11116130708155, -1.0678265669046618, 147.34646739130434, 1.264945652173913]
  e=[-7.347103241623678, -14.93483561943059, 4.885847112946526, 21.54749924818453, 
     -5.148057565540175, 8.136928580553082, -27.90674208255712, 21.23889269084667, 
     292.95889431249236, -0.20306256630643107, -293.11684782608694, -0.11141304347826086]
\end{verbatim}


\begin{thebibliography}{}

\bibitem{fehlberg}
E.~Fehlberg, Klassische Runge-Kutta-Formeln vierter und niedrigerer Ordnung mit Schrittweiten-Kontrolle und ihre Anwendung auf Wärmeleitungsprobleme, 
Computing 6, 61–71 (1970). https://doi.org/10.1007/BF02241732

\bibitem{hairer} E.~Hairer and G.~Wanner, Solving Ordinary Differential Equations II, Stiff and differential algebraic Problems, 
(2nd ed.), Springer-Verlag, Berlin Heidelberg (1996) 

\bibitem{jax2} T.~Jax, A rooted-tree based derivation of ROW-type methods with arbitrary jacobian entries for solving index-one DAEs, 
Dissertation, University Wuppertal (2019)

\bibitem{jaxstei}
T.~Jax, G.~Steinebach,
Generalized ROW-Type Methods for  for solving semi-explicit DAEs of index-1,
Journal of Computational and Applied Mathematics, Vol.316, S. 213-228 (2017)

\bibitem{lang} J.~Lang, Rosenbrock-Wanner Methods: Construction and Mission, 
In: Jax T., Bartel A., Ehrhardt M., Günther M., Steinebach G. (eds) Rosenbrock—Wanner–Type Methods. Mathematics Online First Collections,
Springer, Cham., 1-17, https://doi.org/10.1007/978-3-030-76810-2\_2 (2021)

\bibitem{lang2} J.~Lang, J.G.~Verwer, ROS3P—An Accurate Third-Order Rosenbrock Solver Designed for
  Parabolic Problems, J. BIT Numerical Mathematics (2001) 41: 731. doi:10.1023/A:1021900219772

\bibitem{oro} A.~Ostermann, M.~Roche, Rosenbrock methods for partial differential equations and fractional orders of convergence,
SIAM J. Numer. Anal. 30, 1084-1098 (1993)

\bibitem{julia} C.~Rackauckas, Q.~Nie, Differentialequations.jl--a performant and feature-rich ecosystem for solving differential equations in julia,
Journal of Open Research Software, 5(1), p.15, Ubiquity Press (2017)

\bibitem{rang2}
J.~Rang, The Prothero and Robinson example: Convergence studies for Runge--Kutta and Rosenbrock--Wanner methods,
Applied Numerical Mathematics 108, 37-56 (2016)

\bibitem{renrostei}
P.~Rentrop, M.~Roche, G.~Steinebach, The application of Rosenbrock-Wanner type methods with stepsize control in differential algebraic equations,
Numer. Math. 55, 545-563 (1989) 

\bibitem{roche} M.~Roche, Rosenbrock methods for differential algebraic equations, Numerische Mathematik, 52, 45-63 (1988)

\bibitem{scholz} S.~Scholz, Order barriers for the B-Convergence of ROW Methods, Computing 41, 219-235 (1989)

\bibitem{rodas4p} G.~Steinebach, Order-reduction of ROW-methods for DAEs and method of lines  applications. Preprint-Nr. 1741, FB Mathematik, TH Darmstadt (1995)

\bibitem{rodas4p2} G.~Steinebach, Improvement of Rosenbrock-Wanner Method RODASP, In: Reis T., Grundel S., Schöps S. (eds) 
Progress in Differential-Algebraic Equations II. Differential-Algebraic Equations Forum. Springer, Cham., 165-184,
doi:10.1007/978-3-030-53905-4\_6 (2020)

\bibitem{rodas5p} G.~Steinebach, Construction of Rosenbrock-Wanner method Rodas5P and numerical benchmarks within the Julia Differential Equations package, 
Bit Numer Math 63, 27 (2023). https://doi.org/10.1007/s10543-023-00967-x

\bibitem{steijulia}
G.~Steinebach, Rosenbrock methods within OrdinaryDiffEq.jl - Overview, recent developments and applications -, 
JuliaCon Proceedings, submitted (2024)

\bibitem{tsit5}
Ch.~Tsitouras, Runge–Kutta pairs of order 5(4) satisfying only the first column simplifying assumption,
Computers and Mathematics with Applications 62, 770–775 (2011) 
 
\end{thebibliography}
\end{document}